\documentclass{amsart}
\usepackage{amsmath,amsfonts,amssymb,amsthm,epsfig}

\newtheorem{Theorem}{Theorem}[section]
 
\newtheorem{Lemma}[Theorem]{Lemma}

\newtheorem{Corollary}[Theorem]{Corollary}

\newtheorem{Question}{Question}
\newtheorem{Comment}[Theorem]{Comment}
\newtheorem{Definition}[Theorem]{Definition}

\newcommand{\wt}{\operatorname{wt}}

\newcommand{\bigdot}{\cdot}
\newcommand{\bz}{\Bbb{Z}}
\newcommand{\bc}{\Bbb{C}}
\newcommand{\g}{\mathfrak{g}}

\newcommand{\asl}{\widehat{\text{sl}}}

\newcommand{\slt}{\mbox{sl}_3}
\newcommand{\gli}{\mbox{gl}_\infty}
\newcommand{\hgl}{\widehat{\mbox{gl}}}
\newcommand{\Id}{\text{Id}}

\theoremstyle{definition}

\newcommand{\mathfig}[2]{{\hspace{-3pt}\begin{array}{c}%
  \raisebox{-2.5pt}{\includegraphics[width=#1\textwidth]{#2}}
\end{array}\hspace{-3pt}}}

\begin{document}

\title[Three combinatorial models for $\asl_n$ crystals]{Three combinatorial models for $\asl_n$ crystals, with applications to cylindric plane partitions}

\author{Peter Tingley}
\email{pwtingle@math.berkeley.edu}
\address{UC Berkeley, Department of Mathematics\\ Berkeley, CA}
\thanks{This work was supported by the RTG grant DMS-0354321.}

\begin{abstract}
We define three combinatorial models for $\asl_n$ crystals, parametrized by partitions, configurations of beads on an ``abacus", and cylindric plane partitions, respectively. These are reducible, but we can identify an irreducible subcrystal corresponding to any dominant integral highest weight $\Lambda$. Cylindric plane partitions actually parametrize a basis for $V_\Lambda \otimes F$, where $F$ is the space spanned by partitions. We use this to calculate the partition function for a system of random cylindric plane partitions. We also observe a form of rank level duality. Finally, we use an explicit bijection to relate our work to the Kyoto path model.
\end{abstract}

\date{\today}
\maketitle
\tableofcontents

\section{Introduction}

This work was motivated by the Hayashi realization for crystals of level one $\asl_n$ representations, originally developed by Misra and Miwa \cite{MM:1990} using work of Hayashi \cite{Hayashi:1990} (see also \cite{ariki:2000}, Chapter 10). In that realization, the underlying set of the crystal consists of partitions, and the operators $f_i$ act by adding a box to the associated Young diagram. We wondered if there was a similar realization for representations of arbitrary level $\ell$, where the operators $f_i$ would act by adding an $\ell$-ribbon (see Section \ref{lotsastuff}). It turns out that there is. To prove that our construction works, we need a second model, which is based on the abacus used by James and Kerber (\cite{JK:1981}, Chapter 2.7). The abacus model is reminiscent of  a ``Dirac sea,"  and we think it is interesting in its own right. 

The crystals one obtains from the abacus model are not irreducible. However, one can pick out a ``highest" irreducible sub-crystal, so we do have a model for the crystal of any irreducible integrable highest weight representation of $\asl_n$. There is also a reducible sub-crystal whose underlying set is naturally in bijection with the set of cylindric plane partitions with a given boundary. This is our third combinatorial model. We also define a $\gli$ action on the space spanned by cylindric plane partitions, which commutes with our crystal operators. Furthermore, if we use both actions, the space becomes irreducible; in fact, the set of cylindric plane partitions with a given boundary forms a basis for $V_\Lambda \otimes F$, where $V_\Lambda$ is an irreducible representation of some $\asl_n$, and $F$ is the space spanned by all partitions.

In this picture, the weight $|\pi|$ of a cylindric plane partition $\pi$ is easily calculated from the principally graded weight of the corresponding element of $V_\Lambda$ and the size $|\lambda|$ of the corresponding partition $\lambda$. This allows us to calculate the partition function for the system of random cylindric plane partitions studied by Borodin in \cite{Borodin:2006}. Our answer looks quite different from the formula given by Borodin, but we can directly show that they agree. This gives a new link between Borodin's work and the representation theory of $\asl_n$.

There is a symmetry in our model which allows us to consider a given cylindric plane partition as an element of either $V_\Lambda \otimes F$ or $V_{\Lambda'} \otimes F$, where $\Lambda$ is some level $\ell$ highest weight for $\asl_n$, and $\Lambda'$ is a level $n$ highest weight for $\asl_\ell$ determined by $\Lambda$. This implies an identity of $q$-characters, which is our rank-level duality.  Our result is similar to a duality discovered by Frenkel \cite{IFrenkel:1982}. We show how these are related and obtain a new proof of Frenkel's result.

We then relate our work to the Kyoto path model developed by Kashiwara et. al. in \cite{KKMMNN1} and \cite{KKMMNN2} (see \cite{Hong&Kang:2000} for a more recent explanation). We do this by exhibiting an explicit crystal isomorphism between the highest irreducible component of the abacus model and the Kyoto path model for a particular perfect crystal and ground state path.

We finish with some questions. Most notably, it is natural to ask if our crystal structures can be lifted to get representations of $U_q(\asl_n)$. We believe that there should be such a lifting for the space spanned by cylindric plane partitions. This should be similar to the $q$-deformed Fock space studied by  Kashiwara, Miwa, Petersen and Yung in \cite{KMPY:1996}.

Before beginning, we would like to mention a 1991 paper by Jimbo, Misra, Miwa and Okado \cite{JMMO:1991} which contains some results relevant to the present work. In particular, they present a similar realization for the crystal of any irreducible integrable highest weight representation of $\asl_n$.

\subsection{Acknowledgments}
I would like to thank Mark Haiman, Tony Chiang, Brian Rothbock, Alex Woo, Sami Assaf, and everyone else who attended Mark's seminar in 2003-2004. This paper could never have happened without all of your input. I would also like to thank Alexander Braverman, Alejandra Premat, Anne Schilling and Monica Vazirani for useful discussions. Finally, I  would like to thank my advisor Nicolai Reshetikhin for his patience and support.

\subsection{Change log}

The purpose of this section is to record significant changes to this work since publication (v2). 

\begin{description}
\item[Aug  26, 2008] Section 4.2 was corrected. With Definition 4.9 as stated in the published paper, the crystal structure on cylindric plane partitions described in Section 4.2 and the caption to Figure 13 was incorrect. Note that the crystal structure on descending abacus configurations described in Section 3 was correct. However we incorrectly translated this into the language of cylindric plane partitions. Section 4.2 is independent of the rest of the paper, so the rest of our results remain true as originally stated.
\end{description}

\section{Background}

In this section we review some tools we will need. We only include those results most important for the present work, and refer the reader to other sources for more details.

\subsection{The abacus}  \label{lotsastuff}

Here we explain the abacus used by James and Kerber in \cite{JK:1981}.
We start by defining a bijection between partitions and rows of beads. This is essentially the correspondence between partitions and semi-infinite wedge products (see for example \cite{Kac:1990}, Chapter 14); the positions of the beads correspond to the factors in the wedge product. As in \cite{JK:1981}, this one row is transformed into several parallel rows of beads on an abacus.

We use the ``Russian" diagram of a partition, shown in Figure \ref{partition}. Place a bead on the horizontal axis under each down-sloping segment on the edge of the diagram. The corresponding row of beads uniquely defines the partition. Label the horizontal axis so the corners of all boxes are integers, with the vertex at $0$. For a partition $\lambda= (\lambda_1, \lambda_2, \ldots)$, where we define $\lambda_i=0$ for large $i$, the positions of the beads will be $\lambda_i - i + 1/2$. The empty partition corresponds to having beads at all negative positions in $\bz + 1/2$, and none of the positive positions. Adding a box to the partition corresponds to moving a bead one step to the right.

Recall that a {\it ribbon} is a skew partition $\lambda \backslash  \mu$ whose diagram is connected and has at most one box above each position on the horizontal axis (i.e. it is a strip one box wide). Adding an $\ell$-ribbon (that is, a ribbon with $\ell$ boxes) to a partition $\lambda$ moves one bead exactly $\ell$ steps to the right in the corresponding row of beads, possibly jumping over other beads. See Figure \ref{addribbon}. 
\setlength{\unitlength}{0.4cm}

\begin{figure}
\begin{picture}(30,15)
\put(15,0){\line(1,1){15}}
\put(14,1){\line(1,1){12}}
\put(13,2){\line(1,1){11}}
\put(12,3){\line(1,1){10}}
\put(11,4){\line(1,1){9}}
\put(10,5){\line(1,1){7}}
\put(9,6){\line(1,1){5}}
\put(8,7){\line(1,1){3}}
\put(7,8){\line(1,1){3}}
\put(6,9){\line(1,1){3}}
\put(5,10){\line(1,1){1}}

\put(15,0){\line(-1,1){15}}
\put(16,1){\line(-1,1){10}}
\put(17,2){\line(-1,1){9}}
\put(18,3){\line(-1,1){9}}
\put(19,4){\line(-1,1){6}}
\put(20,5){\line(-1,1){6}}
\put(21,6){\line(-1,1){5}}
\put(22,7){\line(-1,1){5}}
\put(23,8){\line(-1,1){4}}
\put(24,9){\line(-1,1){4}}
\put(25,10){\line(-1,1){3}}
\put(26,11){\line(-1,1){2}}
\put(27,12){\line(-1,1){1}}

\put(2.5,0){\line(0,1){12.5}}
\put(3.5,0){\line(0,1){11.5}}
\put(4.5,0){\line(0,1){10.5}}
\put(6.5,0){\line(0,1){10.5}}
\put(9.5,0){\line(0,1){11.5}}
\put(10.5,0){\line(0,1){10.5}}
\put(11.5,0){\line(0,1){9.5}}
\put(14.5,0){\line(0,1){10.5}}
\put(17.5,0){\line(0,1){11.5}}
\put(20.5,0){\line(0,1){12.5}}
\put(22.5,0){\line(0,1){12.5}}
\put(24.5,0){\line(0,1){12.5}}
\put(26.5,0){\line(0,1){12.5}}

\put(2.5,0){\circle*{0.5}}
\put(3.5,0){\circle*{0.5}}
\put(4.5,0){\circle*{0.5}}
\put(6.5,0){\circle*{0.5}}
\put(9.5,0){\circle*{0.5}}
\put(10.5,0){\circle*{0.5}}
\put(11.5,0){\circle*{0.5}}
\put(14.5,0){\circle*{0.5}}
\put(17.5,0){\circle*{0.5}}
\put(20.5,0){\circle*{0.5}}
\put(22.5,0){\circle*{0.5}}
\put(24.5,0){\circle*{0.5}}
\put(26.5,0){\circle*{0.5}}

\put(1.5,0){\circle*{0.5}}
\put(0.5,0){\circle*{0.5}}
\put(-0.5,0){\circle*{0.5}}

\put(-2.2,-0.05){\ldots}
\end{picture}
\caption{The bead string corresponding to the partition  (12,11,10,9,7,5,3,3,3,1). We fix the origin to be the vertex of the partition, so the beads are in positions
$  \ldots, -12.5, -11.5, -10.5, -8.5, -5.5, -4.5, -3.5, -0.5, 2.5, 5.5, 7.5, 9.5, 11.5.$  \label{partition}}
\end{figure}

\begin{figure}
\begin{picture}(30,15)

\put(15,0){\line(1,1){15}}
\put(14,1){\line(1,1){12}}
\put(13,2){\line(1,1){11}}
\put(12,3){\line(1,1){10}}
\put(11,4){\line(1,1){9}}
\put(10,5){\line(1,1){7}}
\put(9,6){\line(1,1){5}}
\put(8,7){\line(1,1){3}}
\put(7,8){\line(1,1){3}}
\put(6,9){\line(1,1){3}}
\put(5,10){\line(1,1){1}}

\put(15,0){\line(-1,1){15}}
\put(16,1){\line(-1,1){10}}
\put(17,2){\line(-1,1){9}}
\put(18,3){\line(-1,1){9}}
\put(19,4){\line(-1,1){6}}
\put(20,5){\line(-1,1){6}}
\put(21,6){\line(-1,1){5}}
\put(22,7){\line(-1,1){5}}
\put(23,8){\line(-1,1){4}}
\put(24,9){\line(-1,1){4}}
\put(25,10){\line(-1,1){3}}
\put(26,11){\line(-1,1){2}}
\put(27,12){\line(-1,1){1}}

\put(2.5,0){\line(0,1){12.5}}
\put(3.5,0){\line(0,1){11.5}}
\put(4.5,0){\line(0,1){10.5}}
\put(6.5,0){\line(0,1){10.5}}
\put(9.5,0){\line(0,1){11.5}}
\put(10.5,0){\line(0,1){10.5}}
\put(11.5,0){\line(0,1){9.5}}
\put(18.5,0){\line(0,1){12.5}}
\put(17.5,0){\line(0,1){13.5}}
\put(20.5,0){\line(0,1){12.5}}
\put(22.5,0){\line(0,1){12.5}}
\put(24.5,0){\line(0,1){12.5}}
\put(26.5,0){\line(0,1){12.5}}

\put(2.5,0){\circle*{0.5}}
\put(3.5,0){\circle*{0.5}}
\put(4.5,0){\circle*{0.5}}
\put(6.5,0){\circle*{0.5}}
\put(9.5,0){\circle*{0.5}}
\put(10.5,0){\circle*{0.5}}
\put(11.5,0){\circle*{0.5}}
\put(14.5,0){\circle{0.5}}
\put(17.5,0){\circle*{0.5}}
\put(20.5,0){\circle*{0.5}}
\put(22.5,0){\circle*{0.5}}
\put(24.5,0){\circle*{0.5}}
\put(26.5,0){\circle*{0.5}}
\put(18.5,0){\circle*{0.5}}

\put(1.5,0){\circle*{0.5}}
\put(0.5,0){\circle*{0.5}}
\put(-0.5,0){\circle*{0.5}}

\put(-2.2,-0.05){\ldots}

\put(14,11){\line(1,1){3}}
\put(19,12){\line(-1,1){2}}

\end{picture}
\caption{Adding a 4-ribbon corresponds to moving one bead forward 4 positions, possibly jumping over some beads in between. In this example, the bead that was at position $-0.5$ is moved forward four places to position $3.5$. The diagram has changed above positions $0.5, 1.5$ and $2.5$, but adding the ribbon won't change the slope of the diagram there, so the beads stay in the same place.
\label{addribbon}}
\end{figure}
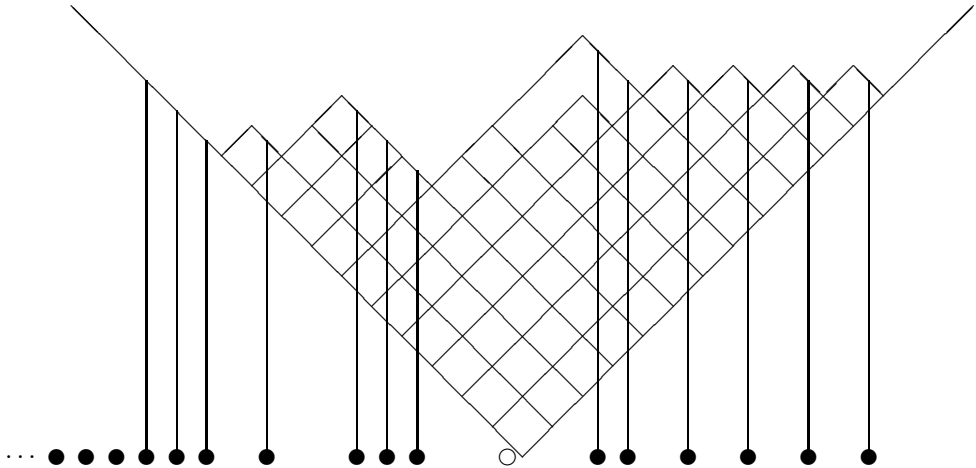

\setlength{\unitlength}{0.5cm}

It is often convenient to work with just the ``bead" picture. In order to avoid confusion, we denote the empty spaces by white beads, and indicate the position of the origin by a line. The example shown in Figure \ref{partition} becomes:

\vspace{-0.22in}

\begin{picture}(20,2)
\label{bead_strand}
\put(5,0){
\begin{picture}(12,0.5)

\put(-5,0){\ldots}
\put(24,0){\ldots}

\put(-3.5,0){\circle*{0.5}}
\put(-2.5,0){\circle*{0.5}}
\put(-1.5,0){\circle*{0.5}}
\put(-0.5,0){\circle*{0.5}}
\put(0.5,0){\circle{0.5}}
\put(1.5,0){\circle*{0.5}}
\put(2.5,0){\circle{0.5}}
\put(3.5,0){\circle{0.5}}
\put(4.5,0){\circle*{0.5}}
\put(5.5,0){\circle*{0.5}}
\put(6.5,0){\circle*{0.5}}
\put(7.5,0){\circle{0.5}}
\put(8.5,0){\circle{0.5}}
\put(9.5,0){\circle*{0.5}}

\put(10,-0.5){\line(0,1){1}}

\put(10.5,0){\circle{0.5}}
\put(11.5,0){\circle{0.5}}
\put(12.5,0){\circle*{0.5}}
\put(13.5,0){\circle{0.5}}
\put(14.5,0){\circle{0.5}}
\put(15.5,0){\circle*{0.5}}
\put(16.5,0){\circle{0.5}}
\put(17.5,0){\circle*{0.5}}
\put(18.5,0){\circle{0.5}}
\put(19.5,0){\circle*{0.5}}
\put(20.5,0){\circle{0.5}}
\put(21.5,0){\circle*{0.5}}
\put(22.5,0){\circle{0.5}}
\put(23.5,0){\circle{0.5}}

\end{picture}

}
\end{picture}
\vspace{0.2in}

\vspace{-0.08in}

\noindent We then put the beads into groups of $\ell$ starting at the origin, as shown below for $\ell=4$:

\vspace{-0.18in}

\begin{picture}(20,2)
\label{grouped_beads}
\put(5,0){
\begin{picture}(12,0.5)

\put(-5,0){\ldots}
\put(24,0){\ldots}

\put(-3.5,0){\circle*{0.5}}
\put(-2.5,0){\circle*{0.5}}

\put(-2,-0.3){\huge(}
\put(-2.4,-0.3){\huge)}

\put(-1.5,0){\circle*{0.5}}
\put(-0.5,0){\circle*{0.5}}
\put(0.5,0){\circle{0.5}}
\put(1.5,0){\circle*{0.5}}

\put(2,-0.3){\huge(}
\put(1.6,-0.3){\huge)}

\put(2.5,0){\circle{0.5}}
\put(3.5,0){\circle{0.5}}
\put(4.5,0){\circle*{0.5}}
\put(5.5,0){\circle*{0.5}}

\put(6,-0.3){\huge(}
\put(5.6,-0.3){\huge)}

\put(6.5,0){\circle*{0.5}}
\put(7.5,0){\circle{0.5}}
\put(8.5,0){\circle{0.5}}
\put(9.5,0){\circle*{0.5}}

\put(10,-0.3){\huge(}
\put(9.6,-0.3){\huge)}

\put(10,-0.5){\line(0,1){1}}

\put(10.5,0){\circle{0.5}}
\put(11.5,0){\circle{0.5}}
\put(12.5,0){\circle*{0.5}}
\put(13.5,0){\circle{0.5}}

\put(14,-0.3){\huge(}
\put(13.6,-0.3){\huge)}

\put(14.5,0){\circle{0.5}}
\put(15.5,0){\circle*{0.5}}
\put(16.5,0){\circle{0.5}}
\put(17.5,0){\circle*{0.5}}

\put(18,-0.3){\huge(}
\put(17.6,-0.3){\huge)}

\put(18.5,0){\circle{0.5}}
\put(19.5,0){\circle*{0.5}}
\put(20.5,0){\circle{0.5}}
\put(21.5,0){\circle*{0.5}}

\put(22,-0.3){\huge(}
\put(21.6,-0.3){\huge)}

\put(22.5,0){\circle{0.5}}
\put(23.5,0){\circle{0.5}}

\end{picture}

}
\end{picture}

\vspace{0.13in}

\noindent Rotating each group 90 degrees counterclockwise, and compressing, we get $\ell$ rows:

\begin{center}
\begin{picture}(20,3.7)

\put(0,0){
\begin{picture}(12,2)

\put(3,0){\ldots}
\put(3,1){\ldots}
\put(3,2){\ldots}
\put(3,3){\ldots}

\put(16,0){\ldots}
\put(16,1){\ldots}
\put(16,2){\ldots}
\put(16,3){\ldots}

\put(4.5,3){\circle*{0.5}}
\put(4.5,2){\circle*{0.5}}
\put(4.5,1){\circle*{0.5}}
\put(4.5,0){\circle*{0.5}}

\put(5.5,3){\circle*{0.5}}
\put(5.5,2){\circle*{0.5}}
\put(5.5,1){\circle*{0.5}}
\put(5.5,0){\circle*{0.5}}

\put(6.5,3){\circle*{0.5}}
\put(6.5,2){\circle*{0.5}}
\put(6.5,1){\circle*{0.5}}
\put(6.5,0){\circle*{0.5}}

\put(7.5,3){\circle*{0.5}}
\put(7.5,2){\circle{0.5}}
\put(7.5,1){\circle*{0.5}}
\put(7.5,0){\circle*{0.5}}

\put(8.5,3){\circle*{0.5}}
\put(8.5,2){\circle*{0.5}}
\put(8.5,1){\circle{0.5}}
\put(8.5,0){\circle{0.5}}

\put(9.5,3){\circle*{0.5}}
\put(9.5,2){\circle{0.5}}
\put(9.5,1){\circle{0.5}}
\put(9.5,0){\circle*{0.5}}

\put(10,-0.5){\line(0,1){4}}

\put(10.5,3){\circle{0.5}}
\put(10.5,2){\circle*{0.5}}
\put(10.5,1){\circle{0.5}}
\put(10.5,0){\circle{0.5}}

\put(11.5,3){\circle*{0.5}}
\put(11.5,2){\circle{0.5}}
\put(11.5,1){\circle*{0.5}}
\put(11.5,0){\circle{0.5}}

\put(12.5,3){\circle*{0.5}}
\put(12.5,2){\circle{0.5}}
\put(12.5,1){\circle*{0.5}}
\put(12.5,0){\circle{0.5}}

\put(13.5,3){\circle{0.5}}
\put(13.5,2){\circle{0.5}}
\put(13.5,1){\circle{0.5}}
\put(13.5,0){\circle{0.5}}

\put(14.5,3){\circle{0.5}}
\put(14.5,2){\circle{0.5}}
\put(14.5,1){\circle{0.5}}
\put(14.5,0){\circle{0.5}}

\put(15.5,3){\circle{0.5}}
\put(15.5,2){\circle{0.5}}
\put(15.5,1){\circle{0.5}}
\put(15.5,0){\circle{0.5}}

\end{picture}

}

\end{picture}

\vspace{0.1in}
\end{center}

\noindent We will call this the level $\ell$ abacus (in this example $\ell=4$).
Adding an $\ell$-ribbon now corresponds to moving one bead forward one position, staying on the same row. For instance, adding the ribbon as in Figure \ref{addribbon} corresponds to moving the third black bead from the right on the top row, which gives:

\begin{picture}(20,4)
\put(5,0){
\begin{picture}(12,2)

\put(3,0){\ldots}
\put(3,1){\ldots}
\put(3,2){\ldots}
\put(3,3){\ldots}

\put(16,0){\ldots}
\put(16,1){\ldots}
\put(16,2){\ldots}
\put(16,3){\ldots}

\put(4.5,3){\circle*{0.5}}
\put(4.5,2){\circle*{0.5}}
\put(4.5,1){\circle*{0.5}}
\put(4.5,0){\circle*{0.5}}

\put(5.5,3){\circle*{0.5}}
\put(5.5,2){\circle*{0.5}}
\put(5.5,1){\circle*{0.5}}
\put(5.5,0){\circle*{0.5}}

\put(6.5,3){\circle*{0.5}}
\put(6.5,2){\circle*{0.5}}
\put(6.5,1){\circle*{0.5}}
\put(6.5,0){\circle*{0.5}}

\put(7.5,3){\circle*{0.5}}
\put(7.5,2){\circle{0.5}}
\put(7.5,1){\circle*{0.5}}
\put(7.5,0){\circle*{0.5}}

\put(8.5,3){\circle*{0.5}}
\put(8.5,2){\circle*{0.5}}
\put(8.5,1){\circle{0.5}}
\put(8.5,0){\circle{0.5}}

\put(9.5,3){\circle{0.5}}
\put(9.5,2){\circle{0.5}}
\put(9.5,1){\circle{0.5}}
\put(9.5,0){\circle*{0.5}}

\put(10,-0.5){\line(0,1){4}}

\put(10.5,3){\circle*{0.5}}
\put(10.5,2){\circle*{0.5}}
\put(10.5,1){\circle{0.5}}
\put(10.5,0){\circle{0.5}}

\put(11.5,3){\circle*{0.5}}
\put(11.5,2){\circle{0.5}}
\put(11.5,1){\circle*{0.5}}
\put(11.5,0){\circle{0.5}}

\put(12.5,3){\circle*{0.5}}
\put(12.5,2){\circle{0.5}}
\put(12.5,1){\circle*{0.5}}
\put(12.5,0){\circle{0.5}}

\put(13.5,3){\circle{0.5}}
\put(13.5,2){\circle{0.5}}
\put(13.5,1){\circle{0.5}}
\put(13.5,0){\circle{0.5}}

\put(14.5,3){\circle{0.5}}
\put(14.5,2){\circle{0.5}}
\put(14.5,1){\circle{0.5}}
\put(14.5,0){\circle{0.5}}

\put(15.5,3){\circle{0.5}}
\put(15.5,2){\circle{0.5}}
\put(15.5,1){\circle{0.5}}
\put(15.5,0){\circle{0.5}}

\end{picture}

}
\end{picture}

\vspace{0.1in}

As explained in (\cite{JK:1981} Chapter 2.7), this model immediately gives some interesting information about the partition: The $\ell$ rows can be interpreted as $\ell$ partitions, using the correspondence between partitions and rows of beads (shifting the origin if necessary). This is known as the $\ell$-quotient of $\lambda$.  We can also consider the partition obtained by pushing the beads on each row as far to the left as they will go, but not changing rows (and only doing finitely many moves). This is the $\ell$-core.

\subsection{Crystals and tensor products}

We use notation as in \cite{Hong&Kang:2000}, and refer the reader to that book for a detailed explanation of crystals. For us, a crystal is a set $B$ associated to a representation $V$ of a symmetrizable Kac-Moody algebra $\g$, along with operators $e_i : B \rightarrow B \cup\{ 0 \}$ and $f_i : B \rightarrow B \cup \{ 0 \}$, which satisfy some conditions. The set $B$ records certain combinatorial data associated to $V$, and the operators $e_i$ and $f_i$ correspond to the Chevalley generators $E_i$ and $F_i$ of $\g$. If $B_\Lambda$ is the crystal of an irreducible, integral highest weight module $V_\Lambda$, then $B_\Lambda$ corresponds to a canonical basis for $V_\Lambda$. That is, we can associate to each $b \in B_\Lambda$ a $v(b) \in V_\Lambda$, such that $\{ v(b) : b \in B \}$ is a basis for $V_\Lambda$. 

Often $B$ will be represented as a colored directed graph whose vertices are the elements of $B$, and we have a $c_i$ colored edge from $b_1$ to $b_2$ if $f_i (b_1) = b_2$. This records all the information about $B$, since $e_i (b_2) = b_1$ if and only if $f_i (b_1)=b_2$. 

The tensor product rule for $\g$ modules leads to a tensor product rule for crystals, which we will now review. We then present an equivalent definition of the tensor product rule using strings of brackets. This fits more closely with our later definitions, and also helps explain why many realizations of crystals (for example, the realization of $\mbox{sl}_n$ crystals using Young tableaux) make use of brackets.

We start by defining three elements in the root lattice of $\g$ associated to each element $b \in B$.
Let $\g$ be a symmetrizable Kac-Moody algebra, with simple roots indexed by $I$. Let $B$ the crystal of an integrable representation $V$ of $\g$. For each $b \in B$ and $i \in I$, define:
\begin{align*}
\varepsilon_i (b) & = \max \{ m : e_i^m (b) \neq 0 \} \\
\varphi_i (b) & = \max \{ m : f_i^m (b) \neq 0 \}.
\end{align*}
These are always finite because $V$ is integrable. 
\begin{Definition}
Let $\Lambda_i$ be the fundamental weight associated to $i \in I$.  For each $b \in B$, define three elements in the weight lattice of $\g$ by:
\begin{enumerate}

\item $\displaystyle \varphi(b):= \sum_{i \in I} \varphi_i (b) \Lambda_i$

\item $\displaystyle \varepsilon(b) := \sum_{i \in I} \varepsilon_i (b) \Lambda_i$

\item $\displaystyle \wt (b):=\varphi(b)-\varepsilon(b)$.

\end{enumerate}
\end{Definition}

\begin{Comment} \label{wtcomment}
It turns out that $\wt(b)$ will always be equal to the weight of the corresponding canonical basis element $v(b)$ (see for example \cite{Hong&Kang:2000}).
\end{Comment}

We now give the tensor product rule for crystals, using conventions from \cite{Hong&Kang:2000}.
If $A$ and $B$ are two crystals, the tensor product $A \otimes B$ is the crystal whose underlying set is $\{ a \otimes b: a \in A, b \in B \}$, with operators $e_i$ and $f_i$ defined by:
\begin{equation} \label{edef} e_i (a \otimes b)=
\begin{cases}
e_i  (a) \otimes b, \quad \text{if}\quad  \varphi_i(a) \geq \varepsilon_i(b)\\
a \otimes e_i  (b),\quad \text{otherwise}
\end{cases} \end{equation}
\begin{equation} \label{fdef} f_i (a \otimes b)=
\begin{cases}
f_i  (a) \otimes b, \quad \text{if} \quad  \varphi_i(a) > \varepsilon_i(b)\\
a \otimes f_i (b),\quad \text{otherwise}.
\end{cases} \end{equation}

\noindent This can be reworded as follows. In this form it is known as the signature rule: 

\begin{Lemma} \label{bdef}
For $b \in B$, let $S_i(b)$ be the string of brackets $) \cdots )( \cdots ($, where the number of  $``)"$ is $\varepsilon_i(b)$ and the number of $``("$ is $\varphi_i(b)$. Then the actions of $e_i$ and $f_i$ on $A \otimes B$ can be calculated as follows: 
\begin{equation*}
e_i (a \otimes b)= 
\begin{cases}
a \otimes e_i(b) \neq 0 \quad \mbox{if the first uncanceled } ``)"  \mbox{ from the right in } S_i(a) S_i(b) \mbox{ comes from } S_i(b) \\
e_i(a) \otimes b \neq 0 \quad \mbox{if the first uncanceled } ``)"  \mbox{ from the right in } S_i(a) S_i(b) \mbox{ comes from } S_i(b) \\
0 \hspace{0.8in} \mbox{if there is no uncanceled } ``)" \mbox{ in } S_i(a)S_i(b)
\end{cases}
\end{equation*}
\begin{equation*}
f_i (a \otimes b)= 
\begin{cases}
f_i(a) \otimes b \neq 0 \quad \mbox{if the first uncanceled } ``("  \mbox{ from the left in } S_i(a) S_i(b) \mbox{ comes from } S_i(a) \\
a \otimes f_i (b) \neq 0 \quad \mbox{if the first uncanceled } ``("  \mbox{ from the left in } S_i(a) S_i(b) \mbox{ comes from } S_i(b) \\
0 \hspace{0.8in} \mbox{if there is no uncanceled } ``(" \mbox{ in } S_i(a)S_i(b)
\end{cases}
\end{equation*}
\end{Lemma}

\begin{proof} This formula for calculating $e_i$ and $f_i$ follows immediately from Equations (\ref{edef}) and (\ref{fdef}). To see that $e_i(a \otimes b) \neq 0$ if there are any uncanceled $``)"$, notice that in this case you always act on a factor that contributes at least one $``)"$, and hence has $\epsilon_i > 0$. By the definition of $\epsilon_i$, $e_i$ does not send this element to $0$. The proof for $f_i$ is similar. \end{proof}
The advantage of Lemma \ref{bdef} over equations (\ref{edef}) and (\ref{fdef}) is that we can easily understand the actions of $e_i$ and $f_i$ on the tensor product of several crystals:

\begin{Corollary} \label{mycrystaldef}
 Let $B_1, \ldots, B_k$ be crystals of integrable representations of $\g$. Let $b_1 \otimes \cdots \otimes b_k \in B_1 \otimes \cdots \otimes B_k$. For each $1 \leq j \leq k$, let $S_i(b_j)$ be the string of brackets $) \cdots )( \cdots ($, where the number of  $``)"$ is $\varepsilon_i(b_j)$ and the number of $``("$ is $\varphi_i(b_j)$. Then:
 \begin{enumerate} 
\item \label{cc1}  
$\displaystyle e_i (b_1 \otimes \cdots \otimes  b_k) =
\begin{cases}
b_1 \otimes  \cdots e_i(b_j)  \cdots \otimes b_k \neq 0 \quad \begin{array}{l}  \mbox{if the first uncanceled } ``)"  \mbox{ from the right in } \\ S_i(b_1) \cdots S_i(b_k) \mbox{ comes from } S_i(b_j) \end{array} \\
0 \hspace{1.6in} \mbox{if there is no uncanceled } ``)" \mbox{ in } S_i(b_1) \cdots S_i(b_k).
\end{cases}$

\item \label{cc2}
$\displaystyle f_i (b_1 \otimes \cdots \otimes  b_k)= 
\begin{cases}
b_1 \otimes  \cdots f_i(b_j)  \cdots \otimes b_k \neq 0 \quad \begin{array}{l}  \mbox{if the first uncanceled } ``("  \mbox{ from the left in } \\ S_i(b_1) \cdots S_i(b_k)  \mbox{ comes from } S_i(b_j) \end{array} \\
0 \hspace{1.6in} \mbox{if there is no uncanceled } ``(" \mbox{ in } S_i(b_1) \cdots S_i(b_k).
\end{cases}$

\item \label{cc3}
$\varepsilon_i(b_1 \otimes \cdots \otimes b_k)$ is the number of uncanceled $``)"$ in $S_i(b_1) \cdots S_i(b_k)$.

\item  \label{cc4}
$\varphi_i (b_1 \otimes \cdots \otimes b_k)$ is the number of uncanceled $``("$ in $S_i(b_1) \cdots S_i(b_k)$.

\end{enumerate}
 \end{Corollary}
 
 \begin{proof}
 Parts (\ref{cc1}) and (\ref{cc2}) follow by iterating Lemma \ref{bdef}. To see part (\ref{cc3}), notice that, if the first uncanceled $``)"$ is in $S_i (b_j)$, then $\varepsilon_i (b_j) \geq 1$ and $e_i$ acts on $b_j$. Hence $e_i$ changes $S_i(b_j)$ by reducing the number of $``)"$ by one, and increasing the number of $``("$ by one. The only affect on $S_i(b_1) \cdots S_i (b_k)$ is that the first uncanceled $``)"$ is changed to $``("$. This reduces the number of uncanceled $``)"$ by one. $e_i$ will send the element to $0$ exactly when there are no uncanceled $``)"$ left. Hence part (\ref{cc3}) follows by the definition of $\varepsilon_i$. Part (\ref{cc4}) is similar. \end{proof}

\section{Crystal structures}  \label{crystal_def}

In this section, we define two families of crystal structures for $\asl_n$. In the first, the vertices are partitions. In the second, the vertices are configurations of beads on an abacus. In each case the family is indexed by a positive integer $\ell$. We will refer to $\ell$ as the level of the crystal, since, using the results of this section, it is straightforward to see our level $\ell$ crystal structure decomposes into a union of crystals corresponding to level $\ell$ irreducible representations of $\asl_n$. We will see the crystal graph of every irreducible $\asl_n$ representation appear as an easily identifiable subcrystal of the abacus model. We do not see every irreducible representation using the partition model as we define it, although we can get the others using a simple change of conventions (just shift the coloring in Figure \ref{ribboncrystal}). There is a slight subtlety in the case of $\asl_2$, since the proof of Theorem \ref{crystal_proof_th} fails. Theorem \ref{hipart} and the results from Section \ref{cylindric_partitions_and_abacus} do hold in this case, as we show in Section \ref{like_kashiwara} using the Kyoto path model. However, we cannot prove the full strength of Theorem \ref{crystal_proof_th} for $n =2$. 

\subsection{Crystal structure on partitions} \label{crystal_for_ribbons}
We now define level $\ell$ crystal operators for $\asl_n$, acting on the set of partitions (at least when $n \geq 3$). Figure \ref{ribboncrystal} illustrates the definition for $n=3$ and $\ell=4$: Color the boxes of a partition with $n$ colors $c_0, c_1, \ldots c_{n-1}$, where all boxes above position $k$ on the horizontal axis are colored $c_s$ for $\displaystyle s \equiv  \huge{ \lfloor} k /l \huge{\rfloor}$ modulo $n$. To act by $f_i$, place a $``("$ above the horizontal position $k$ if boxes in that position are colored $c_i$, and  you  can add an $\ell$-ribbon whose rightmost box is above $k$. Similarly, we put a $``)"$ above each position $k$ where boxes are colored $c_i$ and you can remove an $\ell$-ribbon whose rightmost box is above $k$. $f_i$ acts by adding an $\ell$-ribbon whose rightmost box is below the first uncanceled $``("$ from the left, if possible, and sending that partition to $0$ if there is no uncanceled $``("$. Similarly, $e_i$ acts by removing an $\ell$-ribbon whose rightmost box is below the first uncanceled $``)"$ from the right, if possible, and sending the partition to $0$ otherwise.

Theorem \ref{crystal_proof_th} below will imply that, for $n \geq 3$, these operators do in fact endow the set of partitions with an $\asl_n$ crystal structure. 

\subsection{Crystal structure on abacus configurations}
\label{abcryst}

Since we can identify partitions with certain configurations of beads on the abacus (abacus configurations), the operators $e_i$ and $f_i$ defined above give operators on the set of abacus configurations coming from partitions. We know from Section \ref{lotsastuff} that adding (removing) an $\ell$-ribbon to a partition corresponds to moving one bead forward (backwards) one position in the corresponding $\ell$-strand abacus. This allows us to translate the operators $e_i$ and $f_i$ to the abacus model. We actually get a definition of operators $e_i$ and $f_i$ on the set of all abacus configurations, regardless of whether or not they actually come from partitions (see Figure \ref{abacuscrystal}):

Color the gaps between the columns of beads with $n$ colors, putting  $c_0$ at the origin, and $c_{[ i \mod n ]}$ in the $i^{th}$ gap, counting left to right. We will include the colors in the diagrams by writing in the $c_i$'s below the corresponding gap. The operators $e_i$ and $f_i$ are then calculated as follows: Put a $``("$ every time a bead could move to the right across color $c_i$, and a $``)"$ every time a bead could move to the left across $c_i$. The brackets are ordered moving up each $c_i$ colored gap in turn from left to right. We group all the brackets corresponding to the same gap above that gap. $f_i$ moves the bead corresponding to the first uncanceled ``(" from the left one place forward, if possible, and sends that element to 0 otherwise. Similarly, $e_i$ moves the bead corresponding the the first uncanceled ``)" from the right one space backwards, if possible, and sends the element to $0$ otherwise. 

\setlength{\unitlength}{0.5cm}
\begin{figure}
\begin{picture}(30,17.5)

\put(4.7,15.5){\huge{(}}
\put(14.7,15.5){\huge{(}}

\put(16.7,15.5){\huge{)}}
\put(17.7,15.5){\huge{)}}
\put(27.7,15.5){\huge{(}}
\put(29.7,15.5){\huge{(}}

\put(5, 15){\line(0,-1){4}}
\put(15, 15){\line(0,-1){2}}
\put(17, 15){\line(0,-1){2}}
\put(18, 15){\line(0,-1){3}}
\put(28, 15){\line(0,-1){1}}

\put(29,14){\line(-1,1){1}}
\put(28,15){\line(-1,-1){1}}
\put(27,14){\line(-1,1){1}}
\put(26,15){\line(-1,-1){2}}
\put(-0.3,13.9){$c_2$}
\put(0.7,12.9){$c_2$}
\put(1.7,11.9){$c_2$}
\put(2.7,10.9){$c_0$}
\put(3.7,9.9){$c_0$}
\put(4.7,8.9){$c_0$}
\put(5.7,7.9){$c_0$}
\put(6.7,6.9){$c_1$}
\put(7.7,5.9){$c_1$}
\put(8.7,4.9){$c_1$}
\put(9.7,3.9){$c_1$}
\put(10.7,2.9){$c_2$}
\put(11.7,1.9){$c_2$}
\put(12.7,0.9){$c_2$}
\put(13.7,-0.1){$c_2$}

\put(15.7,-0.1){$c_0$}
\put(16.7,0.9){$c_0$}
\put(17.7,1.9){$c_0$}
\put(18.7,2.9){$c_1$}
\put(19.7,3.9){$c_1$}
\put(20.7,4.9){$c_1$}
\put(21.7,5.9){$c_1$}
\put(22.7,6.9){$c_2$}
\put(23.7,7.9){$c_2$}
\put(24.7,8.9){$c_2$}
\put(25.7,9.9){$c_2$}
\put(26.7,10.9){$c_0$}
\put(27.7,11.9){$c_0$}
\put(28.7,12.9){$c_0$}
\put(29.7,13.9){$c_0$}

\put(5.7,9.9){$c_0$}
\put(6.7,8.9){$c_1$}
\put(7.7,7.9){$c_1$}
\put(8.7,6.9){$c_1$}
\put(9.7,5.9){$c_1$}
\put(10.7,4.9){$c_2$}
\put(11.7,3.9){$c_2$}
\put(12.7,2.9){$c_2$}
\put(13.7,1.9){$c_2$}
\put(14.7,0.9){$c_0$}
\put(15.7,1.9){$c_0$}
\put(16.7,2.9){$c_0$}
\put(17.7,3.9){$c_0$}
\put(18.7,4.9){$c_1$}
\put(19.7,5.9){$c_1$}
\put(20.7,6.9){$c_1$}
\put(21.7,7.9){$c_1$}
\put(22.7,8.9){$c_2$}
\put(23.7,9.9){$c_2$}
\put(24.7,10.9){$c_2$}
\put(25.7,11.9){$c_2$}

\put(7.7,9.9){$c_1$}
\put(8.7,8.9){$c_1$}
\put(9.7,7.9){$c_1$}
\put(10.7,6.9){$c_2$}
\put(11.7,5.9){$c_2$}
\put(12.7,4.9){$c_2$}
\put(13.7,3.9){$c_2$}
\put(14.7,2.9){$c_0$}
\put(15.7,3.9){$c_0$}
\put(16.7,4.9){$c_0$}
\put(17.7,5.9){$c_0$}
\put(18.7,6.9){$c_1$}
\put(19.7,7.9){$c_1$}
\put(20.7,8.9){$c_1$}
\put(21.7,9.9){$c_1$}
\put(22.7,10.9){$c_2$}
\put(23.7,11.9){$c_2$}

\put(8.7,10.9){$c_1$}
\put(9.7,9.9){$c_1$}
\put(10.7,8.9){$c_2$}
\put(11.7,7.9){$c_2$}
\put(12.7,6.9){$c_2$}
\put(13.7,5.9){$c_2$}
\put(14.7,4.9){$c_0$}
\put(15.7,5.9){$c_0$}
\put(16.7,6.9){$c_0$}
\put(17.7,7.9){$c_0$}
\put(18.7,8.9){$c_1$}
\put(19.7,9.9){$c_1$}
\put(20.7,10.9){$c_1$}
\put(21.7,11.9){$c_1$}

\put(12.7,8.9){$c_2$}
\put(13.7,7.9){$c_2$}
\put(14.7,6.9){$c_0$}
\put(15.7,7.9){$c_0$}
\put(16.7,8.9){$c_0$}
\put(17.7,9.9){$c_0$}
\put(18.7,10.9){$c_1$}
\put(19.7,11.9){$c_1$}

\put(13.7,9.9){$c_2$}
\put(14.7,8.9){$c_0$}
\put(15.7,9.9){$c_0$}
\put(16.7,10.9){$c_0$}
\put(17.7,11.9){$c_0$}

\put(14.7,10.9){$c_0$}
\put(15.7,11.9){$c_0$}
\put(16.7,12.9){$c_0$}

\put(15,0){\line(1,1){15}}
\put(14,1){\line(1,1){12}}
\put(13,2){\line(1,1){11}}
\put(12,3){\line(1,1){10}}
\put(11,4){\line(1,1){9}}
\put(10,5){\line(1,1){8}}
\put(9,6){\line(1,1){8}}
\put(8,7){\line(1,1){3}}
\put(7,8){\line(1,1){3}}
\put(6,9){\line(1,1){3}}
\put(5,10){\line(1,1){1}}

\put(15,0){\line(-1,1){15}}
\put(16,1){\line(-1,1){10}}
\put(17,2){\line(-1,1){9}}
\put(18,3){\line(-1,1){9}}
\put(19,4){\line(-1,1){6}}
\put(20,5){\line(-1,1){6}}
\put(21,6){\line(-1,1){6}}
\put(22,7){\line(-1,1){6}}
\put(23,8){\line(-1,1){6}}
\put(24,9){\line(-1,1){4}}
\put(25,10){\line(-1,1){3}}
\put(26,11){\line(-1,1){2}}
\put(27,12){\line(-1,1){1}}

\end{picture}
\caption{The calculation of $f_0$ acting on the partition shown in Figure \ref{addribbon} for $\ell=4$. \label{ribboncrystal}}
\end{figure}

\begin{figure}
\vspace{0.2in}

\begin{picture}(20,4)

\put(3,0){\ldots}
\put(3,1){\ldots}
\put(3,2){\ldots}
\put(3,3){\ldots}

\put(16,0){\ldots}
\put(16,1){\ldots}
\put(16,2){\ldots}
\put(16,3){\ldots}

\put(4.5,3){\circle*{0.5}}
\put(4.5,2){\circle*{0.5}}
\put(4.5,1){\circle*{0.5}}
\put(4.5,0){\circle*{0.5}}

\put(4.8, -1){$c_1$}

\put(5.5,3){\circle*{0.5}}
\put(5.5,2){\circle*{0.5}}
\put(5.5,1){\circle*{0.5}}
\put(5.5,0){\circle*{0.5}}

\put(5.8, -1){$c_2$}

\put(6.5,3){\circle*{0.5}}
\put(6.5,2){\circle*{0.5}}
\put(6.5,1){\circle*{0.5}}
\put(6.5,0){\circle*{0.5}}

\put(6.8, -1){$c_0$}

\put(6.8, 3.7){\big{(}}

\put(7.5,3){\circle*{0.5}}
\put(7.5,2){\circle{0.5}}
\put(7.5,1){\circle*{0.5}}
\put(7.5,0){\circle*{0.5}}

\put(7.8, -1){$c_1$}

\put(8.5,3){\circle*{0.5}}
\put(8.5,2){\circle*{0.5}}
\put(8.5,1){\circle{0.5}}
\put(8.5,0){\circle{0.5}}

\put(8.8, -1){$c_2$}

\put(9.5,3){\circle{0.5}}
\put(9.5,2){\circle{0.5}}
\put(9.5,1){\circle{0.5}}
\put(9.5,0){\circle*{0.5}}

\put(9.8, -1){$c_0$}

\put(9.3, 3.7){\big{(}}
\put(9.8, 3.7){\big{)}}
\put(10.3, 3.7){\big{)}}

\put(10,-0.5){\line(0,1){4}}

\put(10.5,3){\circle*{0.5}}
\put(10.5,2){\circle*{0.5}}
\put(10.5,1){\circle{0.5}}
\put(10.5,0){\circle{0.5}}

\put(10.8, -1){$c_1$}

\put(11.5,3){\circle*{0.5}}
\put(11.5,2){\circle{0.5}}
\put(11.5,1){\circle*{0.5}}
\put(11.5,0){\circle{0.5}}

\put(11.8, -1){$c_2$}

\put(12.5,3){\circle*{0.5}}
\put(12.5,2){\circle{0.5}}
\put(12.5,1){\circle*{0.5}}
\put(12.5,0){\circle{0.5}}

\put(12.5,0.8){$\longrightarrow$}

\put(12.8, -1){$c_0$}

\put(12.55, 3.7){\big{(}}
\put(13.05, 3.7){\big{(}}

\put(13.5,3){\circle{0.5}}
\put(13.5,2){\circle{0.5}}
\put(13.5,1){\circle{0.5}}
\put(13.5,0){\circle{0.5}}

\put(13.8, -1){$c_1$}

\put(14.5,3){\circle{0.5}}
\put(14.5,2){\circle{0.5}}
\put(14.5,1){\circle{0.5}}
\put(14.5,0){\circle{0.5}}

\put(14.8, -1){$c_2$}

\put(15.5,3){\circle{0.5}}
\put(15.5,2){\circle{0.5}}
\put(15.5,1){\circle{0.5}}
\put(15.5,0){\circle{0.5}}

\end{picture}

\vspace{0.1in}
\caption{The crystal operation shown in Figure \ref{ribboncrystal} can easily be described on the level $4$ abacus. This procedure is clearly well defined for any abacus configuration, regardless of wether or not it corresponds to a partition. \label{abacuscrystal}}
\end{figure}

We will now show that the operators $e_i$ and $f_i$ give the set of abacus configurations the structure of an $\asl_n$ crystal, when $n \geq 3$.  It follows immediately that our level $\ell$ operators on the set of partitions also gives us a crystal structure, since the map sending a partition to the corresponding abacus configuration preserves the operators $e_i$ and $f_i$. 

\begin{Theorem} \label{crystal_proof_th}
Fix $n \geq 3$ and $\ell \geq 1$. Define a colored directed graph $G$ as follows:  The vertices of $G$ are all configurations of beads on an $\ell$-strand abacus, which have finitely many empty positions to the left of the origin, and finitely many full positions to the right of the origin. There is a $c_i$-colored edge from $\psi$ to $\phi$ if and only if $f_i (\psi) = \phi$. Then each connected component of $G$ is the crystal graph of some integrable highest weight representation of $\asl_n$. 
\end{Theorem}

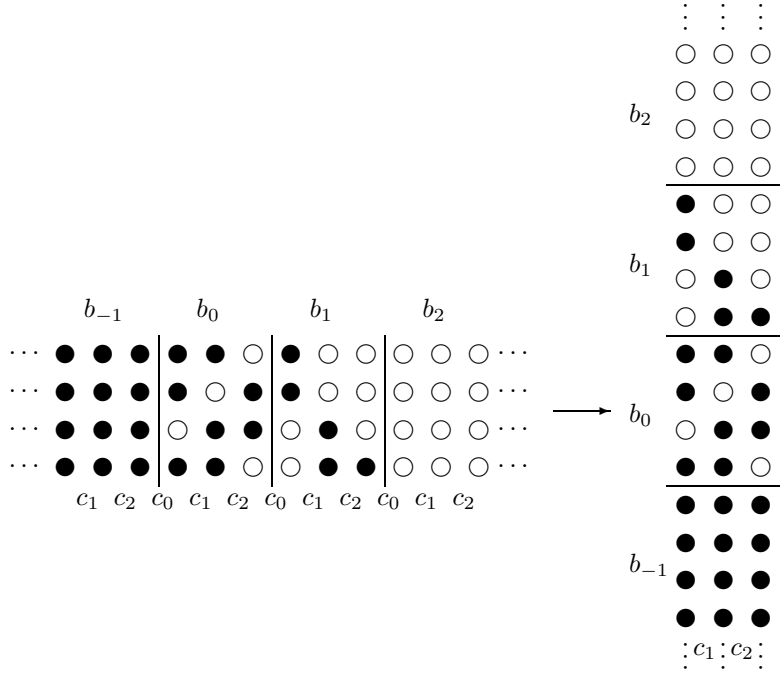
\begin{figure}
\begin{picture}(30,19)

\put(0,-1){\begin{picture}(30,20)

\put(17.5,8.5){\vector(1,0){1.5}}

\put(20.9,18.5){$\cdot$}
\put(20.9,18.8){$\cdot$}
\put(20.9,19.1){$\cdot$}

\put(21.9,18.5){$\cdot$}
\put(21.9,18.8){$\cdot$}
\put(21.9,19.1){$\cdot$}

\put(22.9,18.5){$\cdot$}
\put(22.9,18.8){$\cdot$}
\put(22.9,19.1){$\cdot$}

\put(20.9,1.5){$\cdot$}
\put(20.9,1.8){$\cdot$}
\put(20.9,2.1){$\cdot$}

\put(21.9,1.5){$\cdot$}
\put(21.9,1.8){$\cdot$}
\put(21.9,2.1){$\cdot$}

\put(22.9,1.5){$\cdot$}
\put(22.9,1.8){$\cdot$}
\put(22.9,2.1){$\cdot$}

\put(21, 18){\circle{0.5}}
\put(22, 18){\circle{0.5}}
\put(23, 18){\circle{0.5}}
\put(21, 17){\circle{0.5}}
\put(22, 17){\circle{0.5}}
\put(23, 17){\circle{0.5}}
\put(21, 16){\circle{0.5}}
\put(22, 16){\circle{0.5}}
\put(23, 16){\circle{0.5}}
\put(21, 15){\circle{0.5}}
\put(22, 15){\circle{0.5}}
\put(23, 15){\circle{0.5}}

\put(20.5,14.5){\line(1,0){3}}

\put(21, 14){\circle*{0.5}}
\put(22, 14){\circle{0.5}}
\put(23, 14){\circle{0.5}}
\put(21, 13){\circle*{0.5}}
\put(22, 13){\circle{0.5}}
\put(23, 13){\circle{0.5}}
\put(21, 12){\circle{0.5}}
\put(22, 12){\circle*{0.5}}
\put(23, 12){\circle{0.5}}
\put(21, 11){\circle{0.5}}
\put(22, 11){\circle*{0.5}}
\put(23, 11){\circle*{0.5}}

\put(20.5,10.5){\line(1,0){3}}

\put(21, 10){\circle*{0.5}}
\put(22, 10){\circle*{0.5}}
\put(23, 10){\circle{0.5}}
\put(21, 9){\circle*{0.5}}
\put(22, 9){\circle{0.5}}
\put(23, 9){\circle*{0.5}}
\put(21, 8){\circle{0.5}}
\put(22, 8){\circle*{0.5}}
\put(23, 8){\circle*{0.5}}
\put(21, 7){\circle*{0.5}}
\put(22, 7){\circle*{0.5}}
\put(23, 7){\circle{0.5}}

\put(20.5,6.5){\line(1,0){3}}

\put(21, 6){\circle*{0.5}}
\put(22, 6){\circle*{0.5}}
\put(23, 6){\circle*{0.5}}
\put(21, 5){\circle*{0.5}}
\put(22, 5){\circle*{0.5}}
\put(23, 5){\circle*{0.5}}
\put(21, 4){\circle*{0.5}}
\put(22, 4){\circle*{0.5}}
\put(23, 4){\circle*{0.5}}
\put(21, 3){\circle*{0.5}}
\put(22, 3){\circle*{0.5}}
\put(23, 3){\circle*{0.5}}

\put(21.2,2){$c_1$}
\put(22.2,2){$c_2$}

\put(19.5,4.2){$b_{-1}$}
\put(19.5,8.2){$b_{0}$}
\put(19.5,12.2){$b_{1}$}
\put(19.5,16.2){$b_{2}$}

\put(0,7){\begin{picture}(20,4)
\put(3,0){\ldots}
\put(3,1){\ldots}
\put(3,2){\ldots}
\put(3,3){\ldots}

\put(16,0){\ldots}
\put(16,1){\ldots}
\put(16,2){\ldots}
\put(16,3){\ldots}

\put(4.5,3){\circle*{0.5}}
\put(4.5,2){\circle*{0.5}}
\put(4.5,1){\circle*{0.5}}
\put(4.5,0){\circle*{0.5}}

\put(4.8, -1){$c_1$}

\put(5.5,3){\circle*{0.5}}
\put(5.5,2){\circle*{0.5}}
\put(5.5,1){\circle*{0.5}}
\put(5.5,0){\circle*{0.5}}

\put(5.8, -1){$c_2$}

\put(6.5,3){\circle*{0.5}}
\put(6.5,2){\circle*{0.5}}
\put(6.5,1){\circle*{0.5}}
\put(6.5,0){\circle*{0.5}}

\put(6.8, -1){$c_0$}

\put(7.5,3){\circle*{0.5}}
\put(7.5,2){\circle*{0.5}}
\put(7.5,1){\circle{0.5}}
\put(7.5,0){\circle*{0.5}}

\put(7.8, -1){$c_1$}

\put(8.5,3){\circle*{0.5}}
\put(8.5,2){\circle{0.5}}
\put(8.5,1){\circle*{0.5}}
\put(8.5,0){\circle*{0.5}}

\put(8.8, -1){$c_2$}

\put(9.5,3){\circle{0.5}}
\put(9.5,2){\circle*{0.5}}
\put(9.5,1){\circle*{0.5}}
\put(9.5,0){\circle{0.5}}

\put(9.8, -1){$c_0$}

\put(10,-0.5){\line(0,1){4}}
\put(7,-0.5){\line(0,1){4}}
\put(13,-0.5){\line(0,1){4}}
\put(5,4){$b_{-1}$}
\put(8,4){$b_{0}$}
\put(11,4){$b_{1}$}
\put(14,4){$b_{2}$}

\put(10.5,3){\circle*{0.5}}
\put(10.5,2){\circle*{0.5}}
\put(10.5,1){\circle{0.5}}
\put(10.5,0){\circle{0.5}}

\put(10.8, -1){$c_1$}

\put(11.5,3){\circle{0.5}}
\put(11.5,2){\circle{0.5}}
\put(11.5,1){\circle*{0.5}}
\put(11.5,0){\circle*{0.5}}

\put(11.8, -1){$c_2$}

\put(12.5,3){\circle{0.5}}
\put(12.5,2){\circle{0.5}}
\put(12.5,1){\circle{0.5}}
\put(12.5,0){\circle*{0.5}}

\put(12.8, -1){$c_0$}

\put(13.5,3){\circle{0.5}}
\put(13.5,2){\circle{0.5}}
\put(13.5,1){\circle{0.5}}
\put(13.5,0){\circle{0.5}}

\put(13.8, -1){$c_1$}

\put(14.5,3){\circle{0.5}}
\put(14.5,2){\circle{0.5}}
\put(14.5,1){\circle{0.5}}
\put(14.5,0){\circle{0.5}}

\put(14.8, -1){$c_2$}

\put(15.5,3){\circle{0.5}}
\put(15.5,2){\circle{0.5}}
\put(15.5,1){\circle{0.5}}
\put(15.5,0){\circle{0.5}}

\end{picture}}
\end{picture}}
\end{picture}
\caption{Consider an abacus configuration $\psi$ for $\asl_3$. Cut the abacus along the gaps colored $c_0$ and stack the pieces into a column as shown (this is the step that fails for $\asl_2$). Notice that we can calculate $f_1, f_2, e_1$ and $e_2$ on the column abacus using the same rule as for the row abacus, since we will look at the $c_1$ or $c_2$ gaps in exactly the same order. Each row of the column abacus, considered by itself, generates either a trivial crystal, or the crystal of one of the two fundamental representations of $\slt$. By Corollary \ref{mycrystaldef}, the rule for calculating $e_i$ and $f_i$ on the column abacus is exactly the same as the rule for taking the tensor product of all the non-trivial $\slt$ crystals. So, $\psi$ generates a component of the tensor product of a bunch of crystals of $\slt$ representations, which is itself the crystal of an $\slt$ representation.  \label{crystal_proof_diagram}}
\end{figure}

\begin{proof}
By \cite{KKMMNN1}, Proposition 2.4.4 (see also \cite{Stembridge}), it is sufficient to show that, for each pair $0 \leq i < j<n$, each connected component of the graph obtained by only considering edges of color $c_i$ and $c_j$ is 
\begin{equation*}
\begin{cases} \mbox{An} \hspace{2pt} \mbox{sl}_3 \hspace{2pt} \mbox{crystal if} \hspace{2pt}  |i-j|=1 \hspace{2pt} \mbox{mod}(n) \\  
\mbox{An sl}_2 \times \mbox{sl}_2 \hspace{2pt} \mbox{crystal otherwise.}
\end{cases}
\end{equation*}
Choose some abacus configuration $\psi$, and some $0 \leq i < j \leq n-1$, and consider the connected component containing $\psi$ of the subgraph of $G$ obtained by only considering edges colored $c_i$ and $c_j$. If $|i-j| \neq 1$ $ \mbox{mod} (n)$, then $f_i$ and $f_j$ clearly commute (as do $e_i$ and $e_j$). 
Also, it is clear that if we only consider one color $c_i$, then $G$ is a disjoint union on finite directed lines. This is sufficient to show that the component containing $\psi$ is an $\mbox{sl}_2 \times \mbox{sl}_2$ crystal, as required.

Now consider the case when $|i-j|=1$ $ \mbox{mod} (n)$. The abacus model is clearly symmetric under shifting the colors, so, without loss of generality, we may assume $i=1$ and $j=2$. Also, we need only consider the abacus for $\asl_3$, since if there are columns of beads not bordering a $c_1$ or $c_2$ gap, they can be ignored without affecting $f_1$ or $f_2$. As shown in Figure \ref{crystal_proof_diagram}, each connected component of this crystal is the crystal of an integrable $\slt$ representation, as required. 
\end{proof}

\subsection{More on the abacus model, including the highest irreducible part}
This section is a little technical. We introduce some definitions, then present some results about the structure of the crystals defined in Section \ref{abcryst}. Lemma \ref{structure_lemma} is the main result we need for the applications in Section \ref{cylindric_partitions_and_abacus}. We state and prove this lemma without using the fact that this is an $\asl_n$ crystal, so it holds even for $\asl_2$, when the proof of Theorem \ref{crystal_proof_th} fails. We finish the section by identifying the highest irreducible part, which gives us a realization of the crystal graph for any irreducible integrable representation of $\asl_n$ of level $\ell >0$. We suggest the casual reader look at the statements of Lemma  \ref{structure_lemma} and Theorem \ref{hipart}, and refer to the definitions as needed; the details of the proofs can safely be skipped.

Label the rows of the abacus $0, 1, \ldots, \ell - 1$, where $0$ is at the bottom, 1 is the next row up, and so on. For an abacus configuration $\psi$, we denote the $i^{th}$ row of $\psi$ by $\psi_i$, and the position of the $j^{th}$ bead on that row, counting from the right, by $\psi_i^j$. We will always assume that any row we are considering has only finitely many empty negative positions, and finitely many full positive positions. That is, it can only differ in finitely many places from the row corresponding to the empty partition. 

\begin{Definition}
Let $\psi$ be an abacus configuration. The {\bf compactification} of $\psi$, denoted $\psi_{(0)}$, is the configuration obtained by pushing all the (black) beads to the left, using only finitely many moves, and not changing the row of any bead. For example, Figure \ref{compact1} is the compactification of the configuration in Figure  \ref{good_abacus}.
\end{Definition}

\begin{Definition}
Let $r$ be a row of beads with finitely many negative positions empty and finitely many positive positions full. For each $i \in \bz_{>0}$, let $r^i$ denote the positions of the $i^{th}$ bead in $r$, counting from the right. 
\end{Definition}

\begin{Definition}
Let $r$ and $s$ be two rows of beads. We say $r \leq s$ if $r^i \leq s^i$ for all $i >0$. 
\end{Definition}

\begin{Definition} \label{moveup}
Let $\psi$ be an abacus configuration. Define a row of beads $\psi_i$ for any $i \in \bz$, by letting $\psi_i$ be the $i^{th}$ row of the abacus if $0 \leq i \leq \ell-1$ (counting up and starting with 0), and extending to the rest of $\bz$ using $\psi_{i + \ell}^j := \psi_i^j-n$. That is, $\psi_{i+\ell}$ is $\psi_i$, but shifted $n$ steps to the left.
\end{Definition}

\begin{Definition} \label{descendingdef}
We say a configuration of beads $\psi$ is descending if 
$\psi_0 \geq \psi_1 \geq \ldots \geq \psi_{\ell -1} \geq \psi_{\ell}$. Equivalently, $\psi$ is descending it if $\psi_i \geq \psi_{i+1}$ for all $i \in \bz$.
\end{Definition}

\begin{figure}
\begin{picture}(20,4)

\put(4.5,3){\line(0,-1){2}}
\put(5.5,3){\line(0,-1){2}}
\put(6.5,3){\line(0,-1){2}}

\put(3.5,1){\line(1,-1){1}}
\put(4.5,1){\line(1,-1){1}}
\put(5.5,1){\line(1,-1){1}}
\put(6.5,1){\line(1,-1){1}}

\put(7.5,2){\line(1,-1){1}}
\put(8.5,3){\line(1,-1){1}}
\put(10.5,2){\line(1,-1){1}}
\put(11.5,1){\line(1,-1){1}}

\put(7.5,3){\line(0,-1){1}}
\put(8.5,1){\line(0,-1){1}}
\put(9.5,2){\line(0,-1){1}}
\put(10.5,3){\line(0,-1){1}}

\put(9.5,1){\line(2,-1){2}}

\put(9.33,0.4){$\downarrow$}
\put(9.33,1.4){$\downarrow$}
\put(8.33,2.4){$\downarrow$}
\put(11.33,-0.6){$\downarrow$}
\put(8.33,3.4){$\downarrow$}

\put(3,0){\ldots}
\put(3,1){\ldots}
\put(3,2){\ldots}
\put(3,3){\ldots}

\put(16,0){\ldots}
\put(16,1){\ldots}
\put(16,2){\ldots}
\put(16,3){\ldots}

\put(4.5,3){\circle*{0.5}}
\put(4.5,2){\circle*{0.5}}
\put(4.5,1){\circle*{0.5}}
\put(4.5,0){\circle*{0.5}}

\put(4.8, -1){$c_1$}

\put(5.5,3){\circle*{0.5}}
\put(5.5,2){\circle*{0.5}}
\put(5.5,1){\circle*{0.5}}
\put(5.5,0){\circle*{0.5}}

\put(5.8, -1){$c_2$}

\put(6.5,3){\circle*{0.5}}
\put(6.5,2){\circle*{0.5}}
\put(6.5,1){\circle*{0.5}}
\put(6.5,0){\circle*{0.5}}

\put(6.8, -1){$c_0$}

\put(7.5,3){\circle*{0.5}}
\put(7.5,2){\circle*{0.5}}
\put(7.5,1){\circle{0.5}}
\put(7.5,0){\circle*{0.5}}

\put(7.8, -1){$c_1$}

\put(8.5,3){\circle*{0.5}}
\put(8.5,2){\circle{0.5}}
\put(8.5,1){\circle*{0.5}}
\put(8.5,0){\circle*{0.5}}

\put(8.8, -1){$c_2$}

\put(9.5,3){\circle{0.5}}
\put(9.5,2){\circle*{0.5}}
\put(9.5,1){\circle*{0.5}}
\put(9.5,0){\circle{0.5}}

\put(9.8, -1){$c_0$}

\put(10,-0.5){\line(0,1){4}}

\put(10.5,3){\circle*{0.5}}
\put(10.5,2){\circle*{0.5}}
\put(10.5,1){\circle{0.5}}
\put(10.5,0){\circle{0.5}}

\put(10.8, -1){$c_1$}

\put(11.5,3){\circle{0.5}}
\put(11.5,2){\circle{0.5}}
\put(11.5,1){\circle*{0.5}}
\put(11.5,0){\circle*{0.5}}

\put(11.8, -1){$c_2$}

\put(12.5,3){\circle{0.5}}
\put(12.5,2){\circle{0.5}}
\put(12.5,1){\circle{0.5}}
\put(12.5,0){\circle*{0.5}}

\put(12.8, -1){$c_0$}

\put(13.5,3){\circle{0.5}}
\put(13.5,2){\circle{0.5}}
\put(13.5,1){\circle{0.5}}
\put(13.5,0){\circle{0.5}}

\put(13.8, -1){$c_1$}

\put(14.5,3){\circle{0.5}}
\put(14.5,2){\circle{0.5}}
\put(14.5,1){\circle{0.5}}
\put(14.5,0){\circle{0.5}}

\put(14.8, -1){$c_2$}

\put(15.5,3){\circle{0.5}}
\put(15.5,2){\circle{0.5}}
\put(15.5,1){\circle{0.5}}
\put(15.5,0){\circle{0.5}}

\end{picture}
\vspace{0.1in}
\caption{A descending abacus that is not tight. We have joined the sets $\psi_\cdot^k$ by lines.  $\psi_\cdot^2$ can be shifted as shown. The positions modulo $n$ (here $3$) do not change, the configuration is still decreasing, but they are closer to $\psi_\bigdot^3$. Hence, this configuration is not tight. The operation shown here is $T_2$ as in Definition \ref{dtk}. \label{nottight}}
\end{figure}
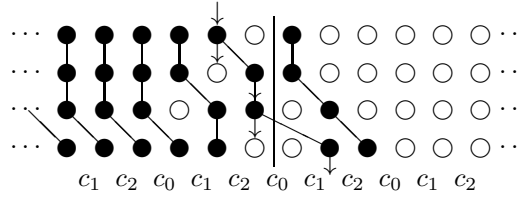

\begin{Comment}
Notice that a descending abacus configuration satisfies $\psi_0^k- n= \psi_\ell^k \leq \psi_{\ell-1}^k$ for all $k$. Thus, if we only draw $\ell$ rows of the abacus, each $\psi_\bigdot^k$ is at most $n+1$ columbs wide (see Figure \ref{nottight}).
\end{Comment}

\begin{Definition} \label{dtk}
The tightening operator $T_k$ is the operator on abacus configuration $\psi$ which shifts the $k^{th}$ bead on each row down one row, if possible. Explicitly:
\begin{equation*}
\begin{cases} 
\begin{cases}
T_k(\psi)_i^j = \psi_i^j \quad \mbox{if} \hspace{3pt} j \neq k \\
T_k(\psi)_i^k= \psi_{i+1}^k
\end{cases}
& \mbox{if} \hspace{3pt} \psi_{i+1}^k > \psi_i^{k+1} \hspace{2pt} \mbox{for all} \hspace{3pt} i \\
T_k (\psi) = 0 & \mbox{otherwise}
\end{cases}
\end{equation*}
See Figure \ref{nottight} for an example. Similarly $T_k^*$ shifts the $k^{th}$ bead on each row up one row:  
\begin{equation*}
\begin{cases} 
\begin{cases}
T_k^*(\psi)_i^j = \psi_i^j \quad \mbox{if} \hspace{3pt} j \neq k \\
T_k^*(\psi)_i^k= \psi_{i-1}^k
\end{cases}
& \mbox{if} \hspace{3pt} \psi_{i-1}^k < \psi_i^{k-1} \hspace{2pt} \mbox{for all} \hspace{3pt} i \\
T_k^* (\psi) = 0 & \mbox{otherwise}
\end{cases}
\end{equation*}
\end{Definition}

\begin{Definition}
We use the notation $\psi_\bigdot^k$ to denote the set consisting of the $k^{th}$ black bead from the right on each row.
\end{Definition}

\begin{Definition} \label{tightdef}
A descending configuration $\psi$ of beads is called {\bf tight} if, for each $k$, the beads $\psi^k_\bigdot$ are positioned as  far to the left as possible, subject to the set
$ \{ \psi^k_\bigdot \mod n \}$ being held fixed, and, for each $i$, $\psi_i^{k+1} < \psi^k_i$.  See Figure \ref{nottight}. This is equivalent to saying $T_k (\psi) =0$ for all $k$, since we only deal with descending configurations, so $T_k$ is the only way to move $\psi_\bigdot^k$ closer to $\psi_\bigdot^{k+1}$.
\end{Definition}

\begin{Lemma} \label{gglemma}
The set of descending abacus configurations (along with 0) is closed under the operators  $e_i$ and $f_i$ defined in Section \ref{abcryst}. Furthermore, the restriction of $e_i$ and $f_i$ to the set of descending abacus configurations can be calculated using the following rule:

For each $k$, let $S_i^k$ be the string of brackets $) \cdots ) ( \cdots ($ where the number of $``("$ is the number of beads in $\psi^k_\bigdot$  in position $i-1/2$ $\mbox{mod}(n)$, and the number of $``)"$ is the number of beads in  $\psi^k_\bigdot$  in position $i+1/2$ $\mbox{mod}(n)$. Let $S_i =  \cdots S_i^3 S_i^2 S_i^1$.

$\bullet$ If the first uncanceled $``)"$ from the right in $S_i$ comes from $S_i^k$, then $e_i$ moves a bead in $\psi_\bigdot^k$ one step to the left. The bead that moves is always the last bead of $\psi_\bigdot^k$ in position  $i+1/2$ $\mbox{mod}(n)$ that you encounter moving up the columns in turn from left to right. If there is no uncanceled $``)"$ in $S_i$, then $e_i$ sends that element to zero.

$\bullet$ If the first uncanceled $``("$ from the left in $S_i$ comes from $S_i^k$, then $f_i$ moves a bead in $\psi_\bigdot^k$ one step to the right. The bead that moves is always the first bead of $\psi_\bigdot^k$ in position  $i-1/2$ $\mbox{mod}(n)$ that you encounter moving up the columns in tun from left to right. If there is no uncanceled $``("$ in $S_i$, then $f_i$ sends that element to zero.
\end{Lemma}

\begin{proof}
Let $R_i$ be the string of brackets used to calculate  $e_i$ and $f_i$ in Section \ref{abcryst}. The descending condition implies:
\begin{enumerate}
\item \label{gg1} All the brackets in $R_i$ coming from beads in $\psi_\bigdot^k$ always come before all the brackets coming from beads in $\psi_\bigdot^{k-1}$.
\item \label{gg2} All the $``)"$ in $R_i$ coming from $\psi_\bigdot^k$ always come before all the $``("$ in $R_i$ coming from $\psi_\bigdot^k$
\end{enumerate}

These facts are both clear if you first do moves as in Figure \ref{cyclic_move} until the first bead of  $\psi_\bigdot^k$ in position $i+1/2$ $\mbox{mod}(n)$ is on the bottom row of the abacus. As argued in the caption, these moves commute with the actions of $e_i$ and $f_i$, and preserves the set of descending abacus configurations.

Let $R_i^k$ be the substring of $R_i$ consisting of brackets coming from  $\psi_\bigdot^k$. 
Then $S_i$ is obtained from $R_i$ by simply adding a string of canceling brackets $( \cdots ( ) \cdots )$ between each $R_i^k$ and $R_i^{k-1}$, where the number of $``("$ and $``)"$ is the number of pairs of touching beads, one in  $\psi_\bigdot^k$ and one in  $\psi_\bigdot^{k-1}$, that are in positions $1-1/2$ and $i+1/2$ $\mbox{mod}(n)$ respectively. Inserting canceling brackets does not change the first uncanceled $``("$.  Hence the first uncanceled $``("$ in $S_i$ will come from a bead in  $\psi_\bigdot^k$ if and only if the first uncanceled $``("$ in $R_i$ comes from a bead in  $\psi_\bigdot^k$. Therefore our new calculation of $f_i$ moves a bead in the right  $\psi_\bigdot^k$. It remains to show that the calculation of $f_i$ using $R_i$ always moves the first bead of  $\psi_\bigdot^k$ is position $i-1/2$ $\mbox{mod} (n)$. But this follows immediately from property (\ref{gg2}) above. The proof for $e_i$ is similar. Hence the new rule agrees with our definition of $e_i$ and $f_i$.

This new rule clearly preserves the set of descending abacus configurations.
\end{proof}

We are now ready to state and prove our main lemma concerning the structure of the operators $e_i$ and $f_i$ acting on abacus configurations:

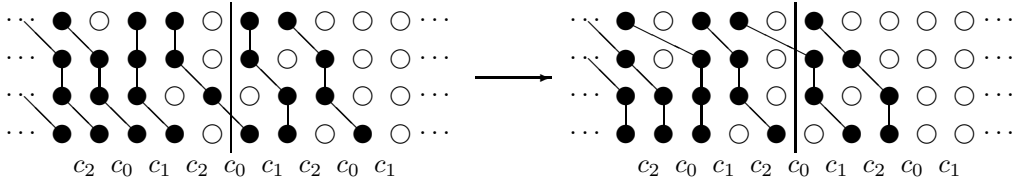
\begin{figure}
\begin{picture}(40,4)
\put(-2,0){\begin{picture}(30,4)

\put(4,0){\ldots}
\put(4,1){\ldots}
\put(4,2){\ldots}
\put(4,3){\ldots}

\put(15,0){\ldots}
\put(15,1){\ldots}
\put(15,2){\ldots}
\put(15,3){\ldots}

\put(4.5,3){\line(1,-1){1}}
\put(5.5,2){\line(0,-1){1}}
\put(5.5,1){\line(1,-1){1}}

\put(5.5,3){\line(1,-1){1}}
\put(6.5,2){\line(0,-1){1}}
\put(6.5,1){\line(1,-1){1}}

\put(7.5,3){\line(0,-1){1}}
\put(7.5,2){\line(0,-1){1}}
\put(7.5,1){\line(1,-1){1}}

\put(8.5,3){\line(0,-1){1}}
\put(8.5,2){\line(1,-1){1}}
\put(9.5,1){\line(1,-1){1}}

\put(10.5,3){\line(0,-1){1}}
\put(10.5,2){\line(1,-1){1}}
\put(11.5,1){\line(0,-1){1}}

\put(11.5,3){\line(1,-1){1}}
\put(12.5,2){\line(0,-1){1}}
\put(12.5,1){\line(1,-1){1}}

\put(4.5,1){\line(1,-1){1}}

\put(5.5,3){\circle*{0.5}}
\put(5.5,2){\circle*{0.5}}
\put(5.5,1){\circle*{0.5}}
\put(5.5,0){\circle*{0.5}}

\put(5.8, -1){$c_2$}

\put(6.5,3){\circle{0.5}}
\put(6.5,2){\circle*{0.5}}
\put(6.5,1){\circle*{0.5}}
\put(6.5,0){\circle*{0.5}}

\put(6.8, -1){$c_0$}

\put(7.5,3){\circle*{0.5}}
\put(7.5,2){\circle*{0.5}}
\put(7.5,1){\circle*{0.5}}
\put(7.5,0){\circle*{0.5}}

\put(7.8, -1){$c_1$}

\put(8.5,3){\circle*{0.5}}
\put(8.5,2){\circle*{0.5}}
\put(8.5,1){\circle{0.5}}
\put(8.5,0){\circle*{0.5}}

\put(8.8, -1){$c_2$}

\put(9.5,3){\circle{0.5}}
\put(9.5,2){\circle{0.5}}
\put(9.5,1){\circle*{0.5}}
\put(9.5,0){\circle{0.5}}

\put(9.8, -1){$c_0$}

\put(10,-0.5){\line(0,1){4}}

\put(10.5,3){\circle*{0.5}}
\put(10.5,2){\circle*{0.5}}
\put(10.5,1){\circle{0.5}}
\put(10.5,0){\circle*{0.5}}

\put(10.8, -1){$c_1$}

\put(11.5,3){\circle*{0.5}}
\put(11.5,2){\circle{0.5}}
\put(11.5,1){\circle*{0.5}}
\put(11.5,0){\circle*{0.5}}

\put(11.8, -1){$c_2$}

\put(12.5,3){\circle{0.5}}
\put(12.5,2){\circle*{0.5}}
\put(12.5,1){\circle*{0.5}}
\put(12.5,0){\circle{0.5}}

\put(12.8, -1){$c_0$}

\put(13.5,3){\circle{0.5}}
\put(13.5,2){\circle{0.5}}
\put(13.5,1){\circle{0.5}}
\put(13.5,0){\circle*{0.5}}

\put(13.8, -1){$c_1$}

\put(14.5,3){\circle{0.5}}
\put(14.5,2){\circle{0.5}}
\put(14.5,1){\circle{0.5}}
\put(14.5,0){\circle{0.5}}

\end{picture}}

\put(14.5, 1.5){\vector(1,0){2}}

\put(13,0){\begin{picture}(30,4)

\put(4,0){\ldots}
\put(4,1){\ldots}
\put(4,2){\ldots}
\put(4,3){\ldots}

\put(15,0){\ldots}
\put(15,1){\ldots}
\put(15,2){\ldots}
\put(15,3){\ldots}

\put(4.5,2){\line(1,-1){1}}
\put(5.5,1){\line(0,-1){1}}

\put(5.5,2){\line(1,-1){1}}
\put(6.5,1){\line(0,-1){1}}

\put(7.5,2){\line(0,-1){1}}
\put(7.5,1){\line(0,-1){1}}

\put(8.5,2){\line(0,-1){1}}
\put(8.5,1){\line(1,-1){1}}

\put(10.5,2){\line(0,-1){1}}
\put(10.5,1){\line(1,-1){1}}

\put(11.5,2){\line(1,-1){1}}
\put(12.5,1){\line(0,-1){1}}

\put(4.5,3){\line(1,-1){1}}
\put(5.5,3){\line(2,-1){2}}
\put(7.5,3){\line(1,-1){1}}
\put(8.5,3){\line(2,-1){2}}
\put(10.5,3){\line(1,-1){1}}

\put(5.5,3){\circle*{0.5}}
\put(5.5,2){\circle*{0.5}}
\put(5.5,1){\circle*{0.5}}
\put(5.5,0){\circle*{0.5}}

\put(5.8, -1){$c_2$}

\put(6.5,3){\circle{0.5}}
\put(6.5,2){\circle{0.5}}
\put(6.5,1){\circle*{0.5}}
\put(6.5,0){\circle*{0.5}}

\put(6.8, -1){$c_0$}

\put(7.5,3){\circle*{0.5}}
\put(7.5,2){\circle*{0.5}}
\put(7.5,1){\circle*{0.5}}
\put(7.5,0){\circle*{0.5}}

\put(7.8, -1){$c_1$}

\put(8.5,3){\circle*{0.5}}
\put(8.5,2){\circle*{0.5}}
\put(8.5,1){\circle*{0.5}}
\put(8.5,0){\circle{0.5}}

\put(8.8, -1){$c_2$}

\put(9.5,3){\circle{0.5}}
\put(9.5,2){\circle{0.5}}
\put(9.5,1){\circle{0.5}}
\put(9.5,0){\circle*{0.5}}

\put(9.8, -1){$c_0$}

\put(10,-0.5){\line(0,1){4}}

\put(10.5,3){\circle*{0.5}}
\put(10.5,2){\circle*{0.5}}
\put(10.5,1){\circle*{0.5}}
\put(10.5,0){\circle{0.5}}

\put(10.8, -1){$c_1$}

\put(11.5,3){\circle{0.5}}
\put(11.5,2){\circle*{0.5}}
\put(11.5,1){\circle{0.5}}
\put(11.5,0){\circle*{0.5}}

\put(11.8, -1){$c_2$}

\put(12.5,3){\circle{0.5}}
\put(12.5,2){\circle{0.5}}
\put(12.5,1){\circle*{0.5}}
\put(12.5,0){\circle*{0.5}}

\put(12.8, -1){$c_0$}

\put(13.5,3){\circle{0.5}}
\put(13.5,2){\circle{0.5}}
\put(13.5,1){\circle{0.5}}
\put(13.5,0){\circle{0.5}}

\put(13.8, -1){$c_1$}

\put(14.5,3){\circle{0.5}}
\put(14.5,2){\circle{0.5}}
\put(14.5,1){\circle{0.5}}
\put(14.5,0){\circle{0.5}}

\end{picture}}
\end{picture}
\vspace{0.02in}
\caption{One can move the bottom row of the abacus to the top row, but shifted back by $n$. This will commute with the operations $f_i$, since when we create the string of brackets to calculate $e_i$ or $f_i$, we will still look at the $c_i$ gaps in exactly the same order. Also, it is clear from Definition \ref{descendingdef} that this preserves the set of descending abacus configurations. \label{cyclic_move}}
\end{figure}

\begin{Lemma} \label{structure_lemma}
Fix $n \geq 2$, and $\ell \geq 1$. Consider the set of $\ell$-strand abacus configurations, colored with $c_0, \ldots c_{n-1}$ as shown in Figure \ref{abacuscrystal}. Let $G$ be the colored directed graph whose vertices are all descending abacus configurations $\psi$ with a given compactification $\psi_{(0)}$, and where there is an edge of color $c_i$ from $\psi$ to $\phi$ if and only if $f_i (\psi)=\phi$. Then:

\begin{enumerate}

\item \label{p2} The operators $e_i$ and $f_i$, restricted to the set of descending abacus configurations, commute with $ T_k$ and $T_k^*$.

\item \label{p3} The sources of $G$ (i.e. vertices that are not the end of any edge) are exactly those configurations that can be obtained from $\psi_{(0)}$ by a series of moves $T_k^*$ for various $k$.
 
\item \label{p3b} The set of tight descending abacus configurations is a connected component of $G$.

\item \label{p4} If we add an edge to $G$ connecting $\psi$ and $\phi$ whenever $T_k (\psi)=\phi$ for some $k$, then $G$ is connected. \end{enumerate}
\end{Lemma}

\begin{proof}

(\ref{p2}): Calculate $f_i$ as in Lemma \ref{gglemma}. Then $T_k$ clearly does not change the string of brackets $S_i$, and hence commute with $f_i$, as long as $T_k (\psi) \neq 0$ and $T_k \circ f_i (\psi) \neq 0$. The only potential problem is if $f_i \circ T_k (\psi) = 0$ but $T_k \circ f_i (\psi) \neq 0$ (or visa versa). 

Assume $f_i \circ T_k (\psi) = 0$ but $T_k \circ f_i (\psi) \neq 0$. Then $f_i$ must move the only bead $b$ of $\psi_\bigdot^k$ that hits $\psi_\bigdot^{k+1}$ when you apply $T_k$. $b$ must be the first bead of $\psi_\bigdot^k$ in position $i-1/2$ $\mbox{mod}(n)$ (moving up the columns and left to right). The last bead of $\psi_\bigdot^k$ in position $i+1/2$ $\mbox{mod}(n)$ must be on the the row below $b$. Since $b$ hits $\psi_\bigdot^{k+1}$ when you apply $T_k$, this must also be the last row of $\psi_\bigdot^{k+1}$ containing a bead in position $i-1/2$ $\mbox{mod}(n)$. By Lemma \ref{gglemma}, for $b$ to move, there must be at least as many beads of $\psi_\bigdot^k$ in position $i+1/2$ $\mbox{mod}(n)$ as beads in $\psi_\bigdot^{k+1}$ in position $i-1/2$ $\mbox{mod}(n)$. The descending condition then implies that the number of beads of $\psi_\bigdot^k$ in position $i+1/2$ $\mbox{mod}(n)$ must be equal to the number of beads of $\psi_\bigdot^{k+1}$ in position $i-1/2$ $\mbox{mod}(n)$. This is illustrated below:
\begin{center}
\begin{picture}(6,5)
\put(2.5,4){\vector(1,0){0.9}}
\put(2.5,3){\line(0,-1){0.7}}
\put(3.5,3){\line(0,-1){0.7}}
\put(2.5,1){\line(0,1){0.7}}
\put(3.5,1){\line(0,1){0.7}}
\put(2.5,4){\circle*{0.5}}
\put(2.5,3){\circle*{0.5}}
\put(2.5,1){\circle*{0.5}}
\put(3.5,3){\circle*{0.5}}
\put(3.5,1){\circle*{0.5}}
\put(2.5,5){\line(0,-1){1}}
\put(1.5,4){\line(1,-1){1}}
\put(2.5,4){\line(1,-1){1}}
\put(2.5,1){\line(1,-1){1}}
\put(3.5,1){\line(1,-1){1}}
\put(2.83,-0.3){$c_i$}
\end{picture}
\end{center}
But then $T_k \circ f_i (\psi)$ is in fact zero, since even after applying $f_i$ the bottom bead of $\psi_\bigdot^k$ in position $i+1/2$ $\mbox{mod}(n)$ cannot be shifted down without hitting $\psi_\bigdot^{k+1}.$ So the problem cannot in fact occur. The other case is similar, as are the cases involving $e_i$ or $T_k^*$.

(\ref{p3}): It is clear from part (\ref{p2}) that  applying operators $T_k^*$ to $\psi_{(0)}$ will always give a source (as long as the configuration is not sent to 0).  So,  let $\psi$ be a source, and we will show that $\psi$ can be tightened to $\psi_{(0)}$ by applying a series of operators $T_k$ for various $k$.  First, for each $k \geq 1$, define:
\begin{align}
 \varphi(\psi_\bigdot^k) := \sum_{j=0}^{\ell-1} \Lambda_i, \mbox{ where } i \equiv \psi_j^k + 1/2 \mbox{ mod} (n).  \\
  \varepsilon(\psi_\bigdot^k) := \sum_{j=0}^{\ell-1} \Lambda_i, \mbox{ where } i \equiv \psi_j^k - 1/2
\mbox{ mod} (n). 
\end{align}
We will use the notation $\varphi_i (\psi)$ (respectively $\varepsilon_i(\psi)$) to mean the coefficient of $\Lambda_i$ in $\varphi (\psi)$ (respectively $\varepsilon(\psi)$). 
For each $k \geq 1$, define $\psi |_k$ to be the abacus configuration obtained by removing the first $k-1$ black beads on each row, counting from the right. When we calculate $f_i$ as in Lemma \ref{gglemma}, all the brackets from $\psi|_k$ always come before the brackets from $\psi_\bigdot^m$, for any $m < k$. Hence, if $\psi$ is a source, then so is $\psi|_k$ for all $k$. We will prove the following statement for each $k \geq 1$:
\begin{align} \label{induct}
 \varphi(\psi|_{k+1}) = \varphi (\psi_\bigdot^{k+1}) = \varepsilon (\psi_\bigdot^{k}). 
\end{align}
$\psi$ can differ in only finitely many places from a compact configuration, so $\psi|_k$ is compact for large enough $k$, and (\ref{induct}) clearly holds for compact configurations. We proceed by induction, assuming (\ref{induct}) holds for some $k \geq 2$ and proving it still holds for $k-1$. 
Since $\varphi(\psi|_{k+1})= \varepsilon (\psi_\bigdot^{k}),$ the rule for calculating $f_i$ implies that $\varphi(\psi|_k)= \varphi(\psi_\bigdot^{k})$. Since $\psi|_{k-1}$ is still a source, we must have $\varepsilon_i(\psi_\bigdot^{k-1}) \leq \varphi_i(\psi_\bigdot^k)$ for each $i$. But $\sum_{i=0}^{n-1} \varepsilon_i(\psi_\bigdot^{k-1}) = \sum_{i=0}^{n-1} \varphi_i(\psi_\bigdot^{k}) = \ell$. Hence we see that  
$ \varphi (\psi_\bigdot^{k}) = \varepsilon (\psi_\bigdot^{k-1})$. So, (\ref{induct}) holds for $k-1$. 

By Definition \ref{tightdef} and the definitions of $\varphi(\psi)$ and $\varepsilon(\psi)$, (\ref{induct}) implies that we can tighten each $\psi_\bigdot^k$ until it is right next to $\psi_\bigdot^{k+1}$. Therefore, $\psi$ can be tightened to a compact configuration, which must be $\psi_{(0)}$ by the definition of $G$. Part (\ref{p3}) follows since $T_k(\psi)=\phi$ if and only if $T_k^*(\phi)=\psi$.

(\ref{p3b}): The graph $G$ is graded by $\bz_{\geq 0}$, where the degree of $\psi$ is the number of times you need to move one bead one step to the left to reach $\psi_{(0)}$. Every connected component to $G$ has at least one vertex $\psi_{\min}$ of minimal degree. $\psi_{\min}$ must be a source, since each $f_i$ is clearly degree 1. Hence each connected component of $G$ contains a source. By part (\ref{p3}), we can tighten any source to get $\psi_{(0)}$. Hence the set of tight descending abacus configurations in $G$ contains exactly one source, namely $\psi_{(0)}$. By Lemma \ref{gglemma} and part (\ref{p2}), the set of tight descending abacus configurations is closed under the operators $f_i$, so it is a complete connected component of $G$. 

(\ref{p4}): This follows immediately from part (\ref{p3}), and the observation from the proof of (\ref{p3b}) that every connected component of $G$ contains a source.
\end{proof}

\begin{Definition} \label{Lambdaa}
Let $\psi_{(0)}$ be a compact, descending abacus configuration. Define a dominant integral weight of $\asl_n$ by
\begin{equation*}
\Lambda(\psi_{(0)}):= \sum_{i=0}^{n-1} m_i \Lambda_i,
\end{equation*}
where $m_i$ is the number of $0 \leq j \leq \ell-1$ such that the last black bead of $\psi_{(0)j}$ is in position $i-1/2$ modulo n. Equivalently, 
\begin{equation*}
m_i = \max \{ m : f_i^m (\psi_{(0)}) \neq 0 \}
\end{equation*}
Note that $\Lambda(\psi_{(0)})$ uniquely determines $\psi_{(0)}$, up to a transformation of the form $\psi_{(0) i}^j \rightarrow \psi_{(0) i+m}^j$ for some $m \in \bz$ (see Definition \ref{moveup}). That is, up to a series of moves as in Figure \ref{cyclic_move}.
\end{Definition}

\begin{Theorem}  \label{hipart} Fix $n \geq 2$ and $\ell \geq 1$, and let $\psi_{(0)}$ be a compact, descending configuration of beads on an $\ell$-strand abacus colored with $c_0, \ldots c_{n-1}$ as shown in Figure \ref{abacuscrystal}. Let $B$ be the colored directed graph whose vertices are all tight, descending abacus configurations with compactification $\psi_{(0)}$, and there is a $c_i$ colored edge from $\psi$ to $\phi$ if $f_i(\psi)=\phi$. Then $B$ is a realization of the crystal graph for the $\asl_n$ representation $V_{\Lambda(\psi_{(0)})}$.
\end{Theorem}

\begin{proof}
Lemma \ref{structure_lemma} part (\ref{p3b}) shows that $B$ is a connected graph. For $n \geq 3$, Theorem \ref{crystal_proof_th} shows that this is in fact the crystal graph of an irreducible $\asl_n$ representation. For $n=2$, it is also the crystal of an irreducible $\asl_2$ representation, but we delay the proof until Section \ref{like_kashiwara}.

It is clear that $\psi_{(0)}$ is the highest weight element. It's weight is $ \sum_{i=0}^{n-1} m_i \Lambda_i$, where $m_i = \max \{ m : f_i^m (\psi_{(0)}) \neq 0 \}$. This is $\Lambda (\psi_{(0)})$ by Definition \ref{Lambdaa}. 
\end{proof}

\begin{figure}
\begin{center}

\begin{picture}(20,3.5)
\put(3,0){\ldots}
\put(3,1){\ldots}
\put(3,2){\ldots}
\put(3,3){\ldots}

\put(16,0){\ldots}
\put(16,1){\ldots}
\put(16,2){\ldots}
\put(16,3){\ldots}

\put(4.5,3){\line(1,-1){1}}
\put(5.5,2){\line(0,-1){1}}
\put(5.5,1){\line(1,-1){1}}

\put(5.5,3){\line(1,-1){1}}
\put(6.5,2){\line(0,-1){1}}
\put(6.5,1){\line(1,-1){1}}
\put(6.5,3){\line(1,-1){1}}
\put(7.5,2){\line(0,-1){1}}
\put(7.5,1){\line(1,-1){1}}
\put(7.5,3){\line(1,-1){1}}
\put(8.5,2){\line(0,-1){1}}
\put(8.5,1){\line(1,-1){1}}
\put(8.5,3){\line(1,-1){1}}
\put(9.5,2){\line(0,-1){1}}
\put(9.5,1){\line(1,-1){1}}
\put(9.5,3){\line(1,-1){1}}
\put(10.5,2){\line(0,-1){1}}
\put(10.5,1){\line(1,-1){1}}
\put(3.5,3){\line(1,-1){1}}
\put(4.5,2){\line(0,-1){1}}
\put(4.5,1){\line(1,-1){1}}
\put(3.5,1){\line(1,-1){1}}

\put(4.5,3){\circle*{0.5}}
\put(4.5,2){\circle*{0.5}}
\put(4.5,1){\circle*{0.5}}
\put(4.5,0){\circle*{0.5}}

\put(4.8, -1){$c_1$}

\put(5.5,3){\circle*{0.5}}
\put(5.5,2){\circle*{0.5}}
\put(5.5,1){\circle*{0.5}}
\put(5.5,0){\circle*{0.5}}

\put(5.8, -1){$c_2$}

\put(6.5,3){\circle*{0.5}}
\put(6.5,2){\circle*{0.5}}
\put(6.5,1){\circle*{0.5}}
\put(6.5,0){\circle*{0.5}}

\put(6.8, -1){$c_0$}

\put(7.5,3){\circle*{0.5}}
\put(7.5,2){\circle*{0.5}}
\put(7.5,1){\circle*{0.5}}
\put(7.5,0){\circle*{0.5}}

\put(7.8, -1){$c_1$}

\put(8.5,3){\circle*{0.5}}
\put(8.5,2){\circle*{0.5}}
\put(8.5,1){\circle*{0.5}}
\put(8.5,0){\circle*{0.5}}

\put(8.8, -1){$c_2$}

\put(9.5,3){\circle*{0.5}}
\put(9.5,2){\circle*{0.5}}
\put(9.5,1){\circle*{0.5}}
\put(9.5,0){\circle*{0.5}}

\put(9.8, -1){$c_0$}

\put(10,-0.5){\line(0,1){4}}

\put(10.5,3){\circle{0.5}}
\put(10.5,2){\circle*{0.5}}
\put(10.5,1){\circle*{0.5}}
\put(10.5,0){\circle*{0.5}}

\put(10.8, -1){$c_1$}

\put(11.5,3){\circle{0.5}}
\put(11.5,2){\circle{0.5}}
\put(11.5,1){\circle{0.5}}
\put(11.5,0){\circle*{0.5}}

\put(11.8, -1){$c_2$}

\put(12.5,3){\circle{0.5}}
\put(12.5,2){\circle{0.5}}
\put(12.5,1){\circle{0.5}}
\put(12.5,0){\circle{0.5}}

\put(12.8, -1){$c_0$}

\put(13.5,3){\circle{0.5}}
\put(13.5,2){\circle{0.5}}
\put(13.5,1){\circle{0.5}}
\put(13.5,0){\circle{0.5}}

\put(13.8, -1){$c_1$}

\put(14.5,3){\circle{0.5}}
\put(14.5,2){\circle{0.5}}
\put(14.5,1){\circle{0.5}}
\put(14.5,0){\circle{0.5}}

\put(14.8, -1){$c_2$}

\put(15.5,3){\circle{0.5}}
\put(15.5,2){\circle{0.5}}
\put(15.5,1){\circle{0.5}}
\put(15.5,0){\circle{0.5}}
\end{picture}
\vspace{0.2in}
\caption{The generator for an irreducible $\asl_3$ crystal of highest weight $\Lambda_0 + 2 \Lambda_1 + \Lambda_2$. As explained in Section \ref{like_kashiwara}, this corresponds to the ground-state path in the Kyoto path model. Each $4$-tuple of beads connected by a line corresponds to an element of $B_4$.  \label{compact1} }
\end{center}
\end{figure}
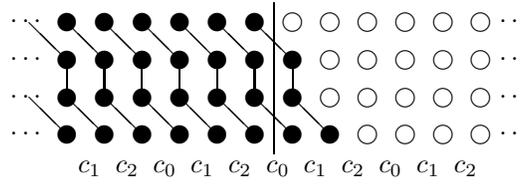

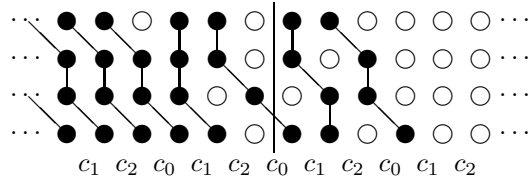
\begin{figure} \vspace{-0.4in}
\begin{center}
\begin{picture}(20,5)
\put(3,0){\ldots}
\put(3,1){\ldots}
\put(3,2){\ldots}
\put(3,3){\ldots}

\put(16,0){\ldots}
\put(16,1){\ldots}
\put(16,2){\ldots}
\put(16,3){\ldots}

\put(3.5,3){\line(1,-1){1}}
\put(4.5,2){\line(0,-1){1}}
\put(4.5,1){\line(1,-1){1}}

\put(4.5,3){\line(1,-1){1}}
\put(5.5,2){\line(0,-1){1}}
\put(5.5,1){\line(1,-1){1}}

\put(5.5,3){\line(1,-1){1}}
\put(6.5,2){\line(0,-1){1}}
\put(6.5,1){\line(1,-1){1}}

\put(7.5,3){\line(0,-1){1}}
\put(7.5,2){\line(0,-1){1}}
\put(7.5,1){\line(1,-1){1}}

\put(8.5,3){\line(0,-1){1}}
\put(8.5,2){\line(1,-1){1}}
\put(9.5,1){\line(1,-1){1}}

\put(10.5,3){\line(0,-1){1}}
\put(10.5,2){\line(1,-1){1}}
\put(11.5,1){\line(0,-1){1}}

\put(11.5,3){\line(1,-1){1}}
\put(12.5,2){\line(0,-1){1}}
\put(12.5,1){\line(1,-1){1}}

\put(3.5,1){\line(1,-1){1}}

\put(4.5,3){\circle*{0.5}}
\put(4.5,2){\circle*{0.5}}
\put(4.5,1){\circle*{0.5}}
\put(4.5,0){\circle*{0.5}}

\put(4.8, -1){$c_1$}

\put(5.5,3){\circle*{0.5}}
\put(5.5,2){\circle*{0.5}}
\put(5.5,1){\circle*{0.5}}
\put(5.5,0){\circle*{0.5}}

\put(5.8, -1){$c_2$}

\put(6.5,3){\circle{0.5}}
\put(6.5,2){\circle*{0.5}}
\put(6.5,1){\circle*{0.5}}
\put(6.5,0){\circle*{0.5}}

\put(6.8, -1){$c_0$}

\put(7.5,3){\circle*{0.5}}
\put(7.5,2){\circle*{0.5}}
\put(7.5,1){\circle*{0.5}}
\put(7.5,0){\circle*{0.5}}

\put(7.8, -1){$c_1$}

\put(8.5,3){\circle*{0.5}}
\put(8.5,2){\circle*{0.5}}
\put(8.5,1){\circle{0.5}}
\put(8.5,0){\circle*{0.5}}

\put(8.8, -1){$c_2$}

\put(9.5,3){\circle{0.5}}
\put(9.5,2){\circle{0.5}}
\put(9.5,1){\circle*{0.5}}
\put(9.5,0){\circle{0.5}}

\put(9.8, -1){$c_0$}

\put(10,-0.5){\line(0,1){4}}

\put(10.5,3){\circle*{0.5}}
\put(10.5,2){\circle*{0.5}}
\put(10.5,1){\circle{0.5}}
\put(10.5,0){\circle*{0.5}}

\put(10.8, -1){$c_1$}

\put(11.5,3){\circle*{0.5}}
\put(11.5,2){\circle{0.5}}
\put(11.5,1){\circle*{0.5}}
\put(11.5,0){\circle*{0.5}}

\put(11.8, -1){$c_2$}

\put(12.5,3){\circle{0.5}}
\put(12.5,2){\circle*{0.5}}
\put(12.5,1){\circle*{0.5}}
\put(12.5,0){\circle{0.5}}

\put(12.8, -1){$c_0$}

\put(13.5,3){\circle{0.5}}
\put(13.5,2){\circle{0.5}}
\put(13.5,1){\circle{0.5}}
\put(13.5,0){\circle*{0.5}}

\put(13.8, -1){$c_1$}

\put(14.5,3){\circle{0.5}}
\put(14.5,2){\circle{0.5}}
\put(14.5,1){\circle{0.5}}
\put(14.5,0){\circle{0.5}}

\put(14.8, -1){$c_2$}

\put(15.5,3){\circle{0.5}}
\put(15.5,2){\circle{0.5}}
\put(15.5,1){\circle{0.5}}
\put(15.5,0){\circle{0.5}}
\end{picture}
\vspace{0.2in}

\end{center}

\caption{This configuration represents an element of the irreducible crystal generated by the compact abacus configuration shown in Figure \ref{compact1}. \label{good_abacus}}

\end{figure}

This implies we can realize any integral irreducible highest weight crystal for $\asl_n$ using the abacus model by choosing the appropriate $\psi_{(0)}$ (see Figure \ref{compact1} for a typical highest weight element $\psi_{(0)}$, and Figure \ref{good_abacus} for a typical element of the highest irreducible part). However, we do not get every irreducible representation from the partition model. To see this, recall that the charge of a row of beads $r$ is the integer $c$ such that, when we push all the black  beads of $r$ to the left, the last black bead is in position $c-1/2$.  Define the charge of an abacus configuration $\psi$ to be the sum of the charges of the rows. The charge of a compact abacus configuration of highest weight $\Lambda$ is then well defined modulo $n$. Only those $\Lambda$ which correspond to compact configurations of charge zero $\mbox{mod} (n)$ can be realized using the crystal structure on partitions, as we define it. However, one can shift the colors in Figure \ref{ribboncrystal} to realize the other irreducible crystals.

\section{Relation to cylindric plane partitions} \label{cylindric_partitions_and_abacus}

We construct a bijection between  the set of descending abacus configurations with a given compactification and the set of cylindric plane partitions on a certain cylinder. This allows us to define a crystal structure on cylindric plane partitions, which is explicitly described in Section \ref{cpp_crystal}. The crystal we obtain is reducible, and carries an action of $\gli$ which commutes with $e_i$ and $f_i$. This action allows us to calculate the partition function for a system of random cylindric plane partitions. We have not seen our formula in the literature, although, as shown in Section \ref{rcpp}, it is equivalent to the formula given by Borodin in \cite{Borodin:2006}. We also observe a form of rank-level duality, which comes from reflecting the cylindric plane partition in a vertical axis. 

\begin{figure} \setlength{\unitlength}{0.3cm}

\begin{center}
\begin{picture}(32,15)

\put(0,-17){ \begin{picture}(32,32)

\put(4,30){$\pi_5$}
\put(6.5,30){$\pi_4$}
\put(9,30){$\pi_3$}
\put(11.5,30){$\pi_2$}
\put(14,30){$\pi_1$}
\put(16.5,30){$\pi_0$}

\put(5,29){\vector(1,-2){0.8}}
\put(7.5,29){\vector(1,-2){0.8}}
\put(10,29){\vector(1,-2){0.8}}
\put(12.5,29){\vector(1,-2){0.8}}
\put(15,29){\vector(1,-2){0.8}}
\put(17.5,29){\vector(1,-2){0.8}}

\put(24,27){$c_2$}
\put(23,29){$c_1$}
\put(22,31){$c_0$}

\put(23.5,26.5){\vector(-2,-1){2}}
\put(22.5,28.5){\vector(-2,-1){2}}
\put(21.5,30.5){\vector(-2,-1){2}}

\put(6.8,24.9){\tiny{1}}
\put(8.8,25.9){\tiny{2}}
\put(11.8,24.9){\tiny{4}}
\put(13.8,25.9){\tiny{5}}
\put(15.8,26.9){\tiny{6}}
\put(18.8,25.9){\tiny{8}}

\put(20.5,25.4){\tiny{9}}
\put(17.5,26.4){\tiny{7}}
\put(10.5,25.4){\tiny{3}}

\linethickness{0.3mm}
\put(-1,32){\line(1,0){29}}
\put(6,32){\line(0,-1){15}}
\put(21,32){\line(0,-1){15}}

\put(6,24.1){\line(2,1){4}}
\put(6,24.05){\line(2,1){4}}
\put(10.05,26){\line(1,-2){1}}
\put(10.1,26){\line(1,-2){1}}
\put(11,24.1){\line(2,1){6}}
\put(11,24.05){\line(2,1){6}}
\put(17.05,27){\line(1,-2){1}}
\put(17.1,27){\line(1,-2){1}}
\put(18,25.1){\line(2,1){2}}
\put(18,25.05){\line(2,1){2}}
\put(20.05,26){\line(1,-2){1}}
\put(20.1,26){\line(1,-2){1}}

\linethickness{0.15mm}

\put(1.25,25.2){6}
\put(2.25,23.2){5}
\put(3.25,21.2){1}

\put(4.25,24.2){7}
\put(5.25,22.2){3}

\put(7.25,23.2){3}
\put(8.25,21.2){3}

\put(9.25,24.2){5}
\put(10.23,22.2){3}

\put(12.25,23.2){4}
\put(13.25,21.2){1}

\put(14.25,24.2){4}
\put(15.25,22.2){4}
\put(16.25,20.2){1}

\put(16.25,25.2){6}
\put(17.25,23.2){5}
\put(18.25,21.2){1}

\put(19.25,24.2){7}
\put(20.25,22.2){3}

\put(22.25,23.2){3}
\put(23.25,21.2){3}

\put(24.25,24.2){5}
\put(25.23,22.2){3}

\put(0,26){\line(1,-2){4}}
\put(2,27){\line(1,-2){4.5}}
\put(5,26){\line(1,-2){4}}
\put(8,25){\line(1,-2){3.5}}
\put(10,26){\line(1,-2){4}}
\put(13,25){\line(1,-2){3.5}}
\put(15,26){\line(1,-2){4}}
\put(17,27){\line(1,-2){4.5}}
\put(20,26){\line(1,-2){4}}
\put(23,25){\line(1,-2){3.5}}
\put(25,26){\line(1,-2){1}}

\put(2,27){\line(-2,-1){3}}
\put(5,26){\line(-2,-1){6}}
\put(10,26){\line(-2,-1){11}}
\put(17,27){\line(-2,-1){18}}
\put(20,26){\line(-2,-1){16}}
\put(25,26){\line(-2,-1){16}}
\put(26,24){\line(-2,-1){12}}
\put(26,21.5){\line(-2,-1){7}}
\put(26,19){\line(-2,-1){2}}

 \end{picture}
}
\end{picture}

\end{center}

\setlength{\unitlength}{0.5cm}
\caption{ A cylindric plane partition with period 9. Here $n=3$ and $\ell=6$, because there are 3 lines on the boundary sloping down and 6 sloping up in each period. $\pi_{ij}$ from Definition \ref{cppd} is the number at the intersection of diagonals $\pi_i$ and $c_j$ (where the empty squares are filled in with $0$). By Theorem \ref{weights} part (\ref{ccpp}), this cylindric plane partition corresponds to a basis element of $V_{2\Lambda_0+3 \Lambda_1 +\Lambda_2} \otimes F$ (a representation of $\asl_3 \oplus \gli$).  We can reflect the picture in a vertical axis, interchanging the $\pi$ diagonals with the $c$ diagonals. In this way the same cylindric plane partition corresponds to a basis element of $V_{\Lambda_0+\Lambda_1+\Lambda_4} \otimes F$ (a representation of $\asl_6 \oplus \gli$). \label{period_nl}}
\end{figure}
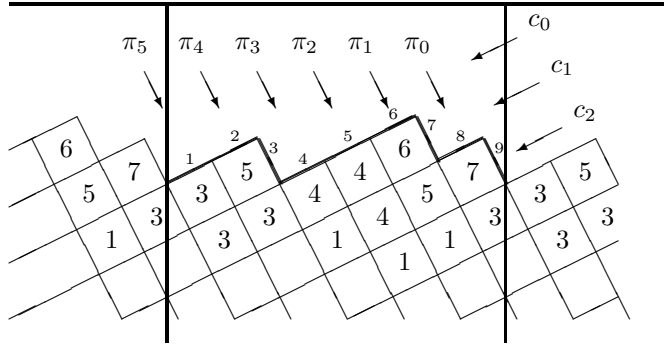

\subsection{A bijection between descending abacus configurations and cylindric plane partitions} \label{bpp}
In this section we define our bijection (Definition \ref{pidef}). Lemma \ref{piisgood} and Theorem \ref{acppb} show that this is a bijection, and relate the compactification of an abacus configuration to the boundary of the corresponding cylindric plane partition. Our definition of a cylindric plane partition (\ref{cppd}) is essentially the same as that used by Borodin in \cite{Borodin:2006}. See Figure \ref{period_nl} for an example. 
  
\begin{Definition} \label{cppd}
A cylindric plane partition of type $(n, \ell)$ is an array of non-negative integers $(\pi_{ij})$, defined for all sufficiently large $i, j \in \bz$, and satisfying:
\begin{enumerate}
\item \label{cppd1} If $\pi_{ij}$ is defined, then so is $\pi_{km}$ whenever $k \geq i$ and $m \geq j$.

\item \label{cppd2} $(\pi_{ij})$ is weakly decreasing in both $i$ and $j$. Furthermore, for all $i$, $\displaystyle \lim_{j \rightarrow \infty} \pi_{ij}=0$, and, for all $j$, $\displaystyle \lim_{i \rightarrow \infty} \pi_{ij}=0$. 

\item \label{cppd3} If $\pi_{ij}$ is defined, then $\pi_{ij} = \pi_{i+\ell, j-n}$.
\end{enumerate}

Notice that condition (\ref{cppd3}) implies that $\pi$ can be represented on a cylinder as in Figure \ref{period_nl}. The period of the cylinder, and the directions of the grid lines, are determined by $n$ and $\ell$. The information of which $\pi_{ij}$ are defined will be called the boundary of the cylinder.
\end{Definition}

\begin{Definition} \label{onemoredef}
A charged partition of charge $k \in \bz$ is a sequence of non-negative integers $(\lambda_{k}  \geq  \lambda_{k+1} \geq \lambda_{k+2} \geq \cdots )$ such that $\lambda_j=0$ for sufficiently large $j$. That is, it is a partition, but with the parts indexed starting at $k$. 
\end{Definition}

\begin{Definition} \label{pidef}
Let $\psi$ be an abacus configuration. For each $i \in \bz$, let $p_i$ be the integer such that the last black bead in the compactification of $\psi_i$ is in position $p_i-1/2$. For $j \geq p_i$, define $\pi_{ij}$ to be the number of black beads to the right of the $j-p_i+1$ st white bead of $\psi_i$, counting from the left.  Define $\pi(\psi):= (\pi_{ij}).$ As shown below (Lemma \ref{piisgood}), $\pi(\psi)$ is a cylindric plane partition. \end{Definition}

\begin{Comment}  \label{alongcomment}
Let $\pi=\pi(\psi)$. It should be clear that each $\pi_{i}:= (\pi_{ip_i}, \pi_{i,p_i+1}, \ldots)$ is a charged partition with charge $c_i$, as defined in Definition \ref{onemoredef}. Also, note that the boundary of $\pi (\psi)$ is determined by the charges of the $\ell$ rows of $\psi$, since $\pi_{ij}$ is well defined exactly when $j$ is at least the charge of $\psi_i$. In particular, the boundary will only depend on the compactification $\psi_{(0)}$ of $\psi$. \end{Comment}

\begin{Lemma} \label{piisgood}
For any descending abacus configuration $\psi$, $\pi(\psi)$ is a cylindric plane partition.
\end{Lemma}

\begin{proof}
Fix a descending abacus configurations $\psi$. 
We first show that part (\ref{cppd2}) of Definition \ref{cppd} holds for $\pi(\psi)$.
As in Section \ref{lotsastuff}, translate each row $\psi_i$ of $\psi$ into a diagram by going down and to the right one step for each black bead, and up and to the right one step for each white bead. Place this diagram at a height so that, far to the right, the diagrams of each $\psi_i$ lie along the same axis, as shown in Figure \ref{cyclic}. Since the $k^{th}$ black bead of $\psi_i$ is always to the right of the $k^{th}$ black bead of $\psi_{i+1}$, one can see that the sequence of diagrams is weakly decreasing by containment (i.e. the diagram of $\psi_{i+1}$ is never above the diagram for $\psi_i$). As in Defnition \ref{pidef}, $\pi_{ij}$ is defined to be the number of black beads to the right of the $j-p_i+1$ st white bead of $\psi_i$. In the diagram, this is the distance of the diagram for $\psi_i$ in the $y$-direction from the interval $[j, j+1]$ on the $x$ axis (see Figure \ref{cyclic}). So $\pi_{ij}$ is weakly decreasing in $i$ because the sequence of diagrams is decreasing. Also $\pi_{ij}$ is weakly decreasing in $j$, since $\pi_i$ is a charged partition. So Definition \ref{cppd} part (\ref{cppd2}) holds. In fact, this argument also shows that Definition \ref{cppd} part (\ref{cppd1}) holds.

Now, by Definition \ref{moveup}, the row $\psi_{i+\ell}$ is just the row $\psi_{i},$ but shifted $n$ steps to the left. By definition \ref{onemoredef}, that means $\pi_{i+\ell}$ will be the same partition as $\pi_i$, but with the charge shifted by $n$. That is, $\pi_{ij}= \pi_{i+\ell, j-n}$, as required in part (\ref{cppd3}) of Definition \ref{cppd}.
\end{proof}

\begin{Definition} \label{Lambdab}
Let $\pi$ be a cylindric plane partition satisfying $\pi_{ij}= \pi_{i+\ell, j-n}$. Label the diagonals $c_0, c_1, c_2, \ldots$, as shown in Figure \ref{cpp}. Let
\begin{equation*}
\Lambda(\pi):= \sum_{i=0}^{n-1} m_i \Lambda_i,
\end{equation*}
where $m_i$ is the number of $0 \leq j \leq \ell-1$ such that the first entry of $\pi_j$ is in diagonal $c_k$ with $k \equiv i$  modulo n. 
\end{Definition}

\begin{Comment}
$\Lambda(\pi)$ is only determined by an (unlabeled) cylindric plane partition $\pi$ up to cyclically relabeling the fundamental weights $\Lambda_i$. This corresponds to the fact that there is a diagram automorphism of $\asl_n$ which cyclically permutes these weights.
\end{Comment}

\setlength{\unitlength}{0.17cm}
\begin{figure}[t]

\begin{center}
\begin{picture}(90,30)
\put(13, 31){$y$}
\put(76,31){$x$}
\put(45,0){\vector(1,1){30}}

\put(45,0){\vector(-1,1){30}}

\put(48,3){\circle*{1}}
\put(51,6){\circle*{1}}
\put(54,9){\circle*{1}}
\put(57,12){\circle*{1}}
\put(60,15){\circle*{1}}
\put(63,18){\circle*{1}}
\put(66,21){\circle*{1}}
\put(69,24){\circle*{1}}

\put(48.5,1.7){$\psi_4 = \psi_{0, (3)}$}
\put(50,4){$\psi_3$}
\put(53,7){$\psi_2$}
\put(54.5,8){$\psi_1$}
\put(57.5,11.5){$\psi_0$}

\put(57.3,12.3){\line(-1,1){3}}
\put(54.3,15.3){\line(-1,-1){3}}
\put(51.3,12.3){\line(-1,1){6}}
\put(45.3,18.3){\line(-1,-1){3}}
\put(42.3,15.3){\vector(-1,1){12}}

\put(54,9){\line(-1,1){6}}
\put(48,15){\line(-1,-1){3}}
\put(45,12){\line(-1,1){3}}
\put(42,15){\line(-1,-1){3}}
\put(39,12){\vector(-1,1){12}}

\put(53.7,8.7){\line(-1,1){3}}
\put(50.7,11.7){\line(-1,-1){3}}
\put(47.7,8.7){\line(-1,1){3}}
\put(44.7,11.7){\line(-1,-1){3}}
\put(41.7,8.7){\vector(-1,1){15}}

\put(50.4,5.4){\line(-1,1){5.7}}
\put(44.7,11.1){\line(-1,-1){3}}
\put(41.7,8.1){\line(-1,1){6}}
\put(35.7,14.1){\line(-1,-1){2.7}}
\put(33,11.4){\vector(-1,1){9}}

\put(48.3,3.3){\line(-1,1){3}}
\put(45.3,6.3){\line(-1,-1){3}}
\put(42.3,3.3){\line(-1,1){6}}
\put(36.3,9.3){\line(-1,-1){3}}
\put(33.3,6.3){\vector(-1,1){12}}
\end{picture}
\caption{
This shows the diagram for $\psi_i$ ($0 \leq i \leq 5$), as described in the proof of Lemma \ref{piisgood}, for the example from Figure \ref{good_abacus}. \label{cyclic}}
\end{center}
\end{figure}
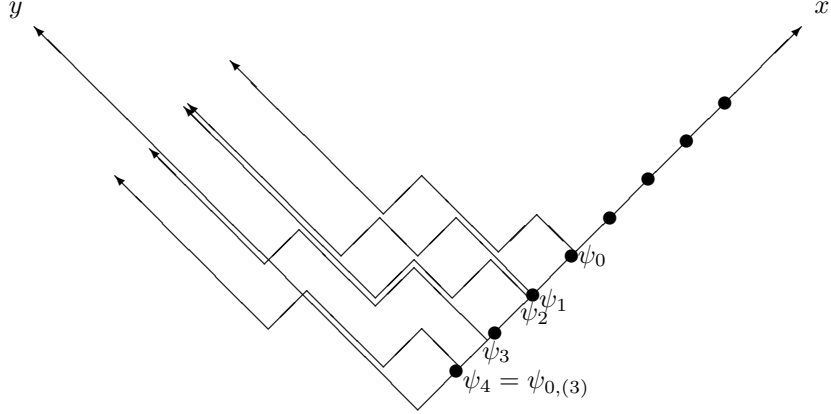

\setlength{\unitlength}{0.3cm}
\begin{figure}
\begin{picture}(40, 16)

\put(13, 10){\line(0,-1){10}}
\put(27, 10){\line(0,-1){10}}

\put(17,15.5){$\pi_0$}
\put(18,15){\vector(1,-1){2}}

\put(15,14.5){$\pi_1$}
\put(16,14){\vector(1,-1){2}}

\put(13,13.5){$\pi_2$}
\put(14,13){\vector(1,-1){2}}

\put(11,12.5){$\pi_3$}
\put(12,12){\vector(1,-1){2}}

\put(9,11.5){$\pi_4$}
\put(10,11){\vector(1,-1){2}}

\put(7,10.5){$\pi_5$}
\put(8,10){\vector(1,-1){2}}

\put(24,14.5){$c_{0}$}
\put(26,13.17){$c_{1}$}
\put(28,11.83){$c_{2}$}
\put(30,10.5){$c_{3}$}
\put(32,9.17){$c_{4}$}

\put(23.5,14){\vector(-2,-1){2}}
\put(25.5,12.67){\vector(-2,-1){2}}
\put(27.5,11.33){\vector(-2,-1){2}}
\put(29.5,10){\vector(-2,-1){2}}
\put(31.5,8.67){\vector(-2,-1){2}}

\put(0,-8){\begin{picture}(40,20)
\put(20,20){}
\put(22,18.67){}
\put(24,17.33){$3$}
\put(26,16){$1$}
\put(28,14.67){$0$}
\put(30,13.33){$0$}
\put(32,12){$0$}

\put(18,19){}
\put(20,17.67){$3$}
\put(22,16.33){$2$}
\put(24,15){$0$}
\put(26,13.67){$0$}
\put(28,12.33){$0$}

\put(16,18){}
\put(18,16.67){$2$}
\put(20,15.33){$1$}
\put(22,14){$0$}
\put(24,12.67){$0$}
\put(26,11.33){$0$}

\put(14,17){$4$}
\put(16,15.67){$2$}
\put(18,14.33){$0$}
\put(20,13){$0$}
\put(22,11.67){$0$}

\put(10,17.33){$3$}
\put(12,16){$1$}
\put(14,14.67){$0$}
\put(16,13.33){$0$}
\put(18,12){$0$}
\put(20,10.67){$0$}

\put(8,16.33){$2$}
\put(10,15){$0$}
\put(12,13.67){$0$}
\put(14,12.33){$0$}
\put(16,11){$0$}

\end{picture} }

\end{picture}

\begin{center}

\caption{$\pi (\psi)$ for the configuration $\psi$ shown in Figures \ref{good_abacus} and \ref{cyclic}. $\pi_{ij}$ is the intersection of the diagonals labeled $\pi_i$ and $c_j$. Notice that $\pi_4$ is just a shift of $\pi_0$. So we can cut on the lines shown, and wrap the diagram around a cylinder, to get a cylindric plane partition. We can also describe Definition \ref{Lambdab}:  $\Lambda (\pi) = \Lambda_0 + 2 \Lambda_1 + \Lambda_2$, since the coefficient of $\Lambda_j$ is the number of $0 \leq i \leq n-1$ for which $\pi_i$ has its first entry at a position $c_k$ with $k \equiv j$ modulo $n$. Note that the choice of labels $c_i$ can be shifted by $c_i \rightarrow c_{i+k}$. This corresponds to rotating the dynkin diagram of $\asl_n$. \label{cpp}}
\end{center}
\end{figure}
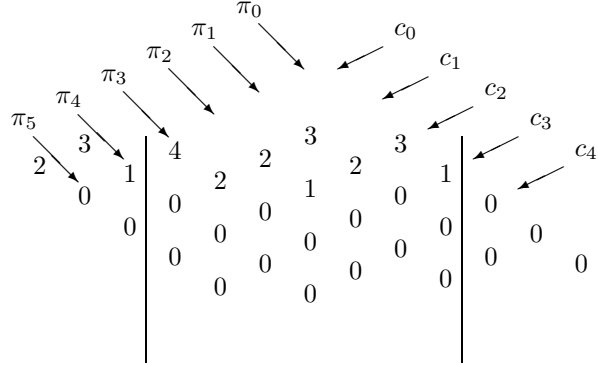
\setlength{\unitlength}{0.5cm}

\begin{Theorem} \label{acppb}
The  map $S: \psi \rightarrow \pi (\psi )$ is a bijection between the set of descending abacus configurations with compactification $\psi_{(0)}$, and the set of cylindric plane partitions with a given boundary. Furthermore, the boundary is determined by $\Lambda(\pi)= \Lambda(\psi_{(0)})$.
\end{Theorem}

\begin{proof}
Lemma \ref{piisgood} shows that $\pi(\psi)$ is always a cylindric plane partition, and, by Comment \ref{alongcomment}, the boundary of $\pi$ depends only on the compactification of $\psi$. Next, notice that, given a cylindric plane partition $\pi$ with the correct boundary, we can construct a diagram as in Figure \ref{cyclic} and then an abacus configuration, simply by reversing the procedure in Definition \ref{pidef}. This is clearly an inverse for $S$, so $S$ is a  bijection. It remains to check that $\Lambda(\pi)= \Lambda(\psi_{(0)})$. The charge of $\pi_i(\psi)$ is both the integer $k$ such that the last black bead of $\psi_i$ is in position $k-1/2$, and the integer $k$ such that the first entry of $\pi_i(\psi)$ is labeled with color $c_k$, so the identity follows from Definitions \ref{Lambdaa} and \ref{Lambdab}.
\end{proof}

\subsection{The crystal structure on cylindric plane partitions} \label{cpp_crystal}

We already have a crystal structure on abacus configurations (for $n \geq 3$) which, by Lemma \ref{gglemma}, preserves the set of descending configurations. Hence Theorem \ref{acppb} implies we have a crystal structure on the set of cylindric plane partitions. At the moment, we need to translate a cylindric plane partition into an abacus configuration to calculate $e_i$ and $f_i$. We now describe how to calculate these operators directly on the cylindric plane partitions. This section is largely independent of the rest of the paper, since later on we find it simpler to work with descending abacus configurations.

\begin{figure}

$$\mathfig{0.5}{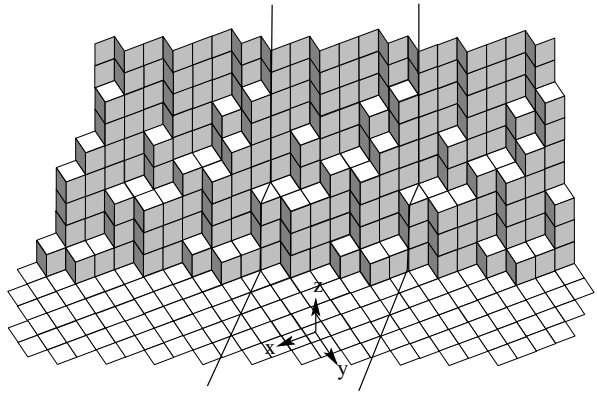}$$

\caption{The three dimensional representation of the cylindric plane partition shown in Figure \ref{period_nl}. The picture is periodic, with one period shown between the dark lines. The first layer of boxes should be colored as in Figure \ref{period_nl}, and higher levels according to the rule that color is constant along $(x_0, y_0+s, z_0+s)$ as $s$ varies. The planes $t(x,y,z)=C$ intersect the ``floor" of the picture in a line which is horizontal in the projection shown, and are angled so that $(x,y,z)$ and $(x,y+1, z+1)$ are always on the same plane. For any $i \in I$, each such plane intersects the center of at most one $c_i$ colored box that could be added or removed  to/from $\pi$ (in each period). The crystal operator $f_i$ acts on the cylindric plane partition by placing a ``('' for each box in $A_i(\pi)$ and a ``)" for each box in $R_i(\pi)$, ordered by $t$ calculated on the coordinates of the center of the box. $f_i$ adds a box corresponding to the first uncanceled ``(", if there is one, and sends the element to $0$ otherwise. Note that adding a box from $R_i(\pi)$ to $\pi$ need not result in a cylindric plane partition. However, it turns out that $f_i(\pi)$ is always a cylindric plane partition.   \label{cpp3}} 
\end{figure}

It will be convenient to view a cylindric plane partition as a 3-dimentional picture, where $\pi_{ij}$ is the height of a pile of boxes placed at position $(i,j)$. See Figure \ref{cpp3}. Each box will be labeled by the coordinates of it's center in the $x$, $y$ and $z$ directions, as shown in Figure \ref{cpp3} (only relative positions matter, so the origin is arbitrary). Each box will also be labeled by a color $c_i$ for some residue $i$ modulo $n$. For the first layer of boxes, this will be determined as in Figure \ref{period_nl}. For higher boxes, one uses the rule that color is constant along lines of the form $\{ (x, y+k, z+k) : k \in \bz \}$. Note that due to the periodicity, $(x,y,z)$ labels the same box as $(x+\ell, y-n, z)$, so the coordinates are only well defined up to this type of transformation. 

For any cylindric plane partition $\pi$ and any residue $i$ modulo $n$, define the following two sets:

\begin{Definition}  $A_i(\pi) $ is the set of all $c_i$ colored boxes that can be added to $\pi$ so that each slice $\pi_i$ (see Figure \ref{cpp}) is still a valid partition.

$R_i(\pi)$ is the set of all $c_i$ colored boxes that can be removed from $\pi$ so that each slice $\pi_i$ is still a valid partition.

{\bf Note:} Adding a box in $A_i(\pi)$ (or removing a box in $R_i(\pi)$) can result in something which is not a cylindric plane partition, since the slices $c_i$ (see Figure \ref{cpp}) may no longer be partitions.
 \end{Definition}

\begin{Definition} \label{defoft}
Define $t(x,y,z)=nx/\ell+y -z.$ Note that $t(x,y,z)= t(x+\ell, y-n, z)$, so $t$ is well defined as a function on boxes in a cylindric plane partition. For a box $n$ as in Figure \ref{cpp3}, define $t(n)$ to be $t$ calculated on the coordinates of the center of $n$.
\end{Definition}

\begin{Lemma} \label{total_order}
Fix a cylindric plane partition $\pi$. Let $n_1, n_2 \in A_i(\pi) \bigcup R_i(\pi)$. Then $t(n_1) \neq t(n_2)$ unless $n_1=n_2$.  
\end{Lemma}

\begin{proof}
Choose $n_0 \in A_i(\pi) \bigcup R_i(\pi)$, and consider the equation $ t(x,y,z)= t(n_0)$. It should be clear from Figure \ref{cpp3} that this plane intersects the center of at most one box in $A_i(\pi) \bigcup R_i(\pi)$ over each period. 
\end{proof}

We are now ready a define our crystal structure. Let $\pi$ be a cylindric plane partition. Define $S_i (\pi)$ to be the string of brackets formed by placing a $``("$ for every box in $A_i(\pi)$, and a $``)"$ for every box in $R_i(\pi)$. These are ordered with the bracket corresponding to $n_1$ coming before the bracket corresponding to $n_2$ if and only if $t(n_1) < t(n_2)$ (this is possible because of Lemma \ref{total_order}). Then $f_i(\pi)$ is the cylindric plane partition obtained by adding the box corresponding to the first uncanceled $``("$, if there is one, and is $0$ otherwise. Simiarly, $e_i(\pi)$ is the cylindric plane partition obtained by removing the box corresponding to the first uncanceled $``)"$, if there is one, and is $0$ otherwise. See Figure \ref{cpp3}.

To see this agrees with the crystal structure coming from descending abacus configurations, one should notice that adding/removing a box $n$ colored $c_i$ with $t(n)=C$ always corresponds to moving a bead across the same gap on the abacus, which we will call $g_i(C)$. One must then check that, for any $c_i$ colored boxes $n_1$ and $n_2$, $g_i(t(n_1))$ comes to the left of $g_i(t(n_2))$ when we calculate the crystal moves on the abacus if and only if $t(n_1)< t(n_2)$. This is straightforward.

\subsection{More structure on cylindric plane partitions}
In this section we give a bijection between descending abacus configurations with compactification $\psi_{(0)}$ and $B_\Lambda \times P$, where $B_\Lambda$ is the crystal associated to an irreducible $\asl_n$ representation $V_\Lambda$, and $P$ is the set of all partitions. It follows that there is also a bijection between cylindric plane partitions with a given boundary and $B_\Lambda \times P$. These bijections are useful because they preserve an appropriate notion of weight.

\begin{Definition}
Let $\psi$ be an abacus configuration. Notice that you can always transform $\psi$ to it's compactification by a series of moves, each of which moves one black bead exactly one step to the left. One can easily see that the number of such moves needed is well defined. This number is the weight of $\psi$, denoted $| \psi |$.
\end{Definition}

\begin{Definition}
If $\pi$ is a cylindric plane partition, the weight of $\pi$, denoted $|\pi |$, is the sum of the entries of $\pi$ over one period. 
\end{Definition}

\begin{Definition} \label{zdef} To each descending abacus configuration $\psi$, assign:
\begin{enumerate}

\item \label{z1} A tight descending abacus configuration $\gamma (\psi)$.

\item \label{z2} A partition $\lambda (\psi)$.

\item \label{z4} A canonical basis vector $v (\psi) \in V_\Lambda$.

\end{enumerate}

\noindent as follows:

(\ref{z1}): Moving from left to right, act by $T_k$ (see definition \ref{dtk}) until $\psi_\bigdot^k$ is tight with respect to $\psi_\bigdot^{k+1}$ for every $k$. The result is $\gamma (\psi)$. 

(\ref{z2}): Construct a row of beads by putting a black bead for each $\psi^i_\bigdot$, and, between the beads corresponding to $\psi^i_\bigdot$ and $\psi^{i+1}_\bigdot$, put a white bead for each time you can apply $T_i$ to $\psi^i_\bigdot$ and stay strictly to the right of $\psi^{i+1}_\bigdot$. Then use the correspondence between rows of beads and partitions to get $\lambda (\psi)$ (shifting to get charge 0 if necessary). 

(\ref{z4}): By Theorem \ref{hipart},  we know that the set of tight descending abacus configurations with compactification $\psi_{(0)}$ are the vertices of the crystal graph for $V_\Lambda$. These in turn correspond to canonical basis vectors in $V_\Lambda$ (see, for example, \cite{Hong&Kang:2000} for a discussion of canonical bases). Let $v(\psi)$ be the canonical basis element associated to $\gamma(\psi)$.
\end{Definition}

\begin{Comment} Since we have a bijection between cylindric plane partitions and abacus configurations, we can also associate this data to a cylindric plane partition $\pi$. In that case, we use the notation $\gamma(\pi), \lambda (\pi)$ and $v(\pi)$. 
\end{Comment}

Let $\gli$ be the Lie algebra of $(\bz+1/2) \times (\bz+1/2)$ matrices with finitely many non-zero entries. As in \cite{Kac:1990}, Chapter 14 (among others), the space $F$ spanned by all partitions is an irreducible $\gli$ module. The generators $E_{p, p+1}$ and $E_{p+1, p}$ of $\gli$ act on a partition as follows: Let $r \in F$, and consider the row of beads corresponding to $r$. 

$\bullet$  If position $p+1$ of $r$ is full and position $p$ is empty then $E_{p, p+1}$ moves the bead in position $p+1$ to position $p$. Otherwise $E_{p, p+1} (r)=0$.

$\bullet$  If position $p$ of $r$ is full and position $p+1$ is empty, then $E_{p+1, p}$ moves the bead in position $p$ to position $p+1$. Otherwise $E_{p+1, p} (r)=0$.

We then get a $\gli$ action on the space of abacus configurations, where $\gli$ acts on $\lambda (\psi)$ without affecting $\gamma(\psi)$. The following theorem will be our key tool in this section:

\begin{Theorem} \label{weights} 
\begin{enumerate}
\item \label{ccommutes} The $\asl_n$ crystal structure on the space of descending abacus configurations with compactification $\psi_{(0)}$ commutes with the $\gli$ action. Furthermore, $ \psi \rightarrow (\gamma (\psi),  \lambda (\psi))$ is a bijection with $B_\Lambda \times P$. Here $\Lambda = \Lambda(\psi_{(0)})$ (see Definition \ref{Lambdaa}), $B_\Lambda$ is the crystal graph of the $\asl_n$ representation $V_\Lambda$ realized according to Theorem \ref{hipart}, and $P$ is the set of all partitions.

\item \label{cabacus} $\{ v(\psi) \otimes \lambda (\psi) \},$ where $\psi$ ranges over all descending abacus configurations with a given compactification $\psi_{(0)}$, is a basis  $V_\Lambda \otimes F$. Here $F$ is the space spanned by all partitions and $\Lambda = \Lambda (\psi_{(0)})$ (see Definition \ref{Lambdaa}). Furthermore, the weight of $\psi$ is the sum of the principle graded weight of $v (\psi)$ and $n |\lambda(\psi)|$.

\item \label{ccpp} $\{ v(\pi) \otimes \lambda(\pi) \}$, where $\pi$ ranges over the set of all cylindric plane partitions $\pi$ with given boundary, is a basis for $V_\Lambda  \otimes F$. Here $F$ is the space of all partitions and $\Lambda = \Lambda (\pi)$ (see Definition \ref{Lambdab}). Furthermore, the weight of such a configuration is the sum of the principle graded weight of $v(\pi)$ and $n |\lambda(\pi)|$.

\end{enumerate}
\end{Theorem}

\begin{proof}
(\ref{ccommutes}): Notice that $\psi$ is uniquely defined by the pair $(\gamma (\psi), \lambda (\psi))$. Lemma \ref{structure_lemma} part (\ref{p2}) implies that $e_i$ and$f_i$ act on $\gamma (\psi)$, ignoring $\lambda (\psi)$. Also, by definition, $\gli$ acts only on $\lambda (\psi)$, ignoring $\gamma (\psi)$. Together, this implies that all $f_i$ commute with the $\gli$ action. Similarly, the $e_i$ commute with the $\gli$ action. Theorem \ref{hipart} says that the $\asl_n$ crystal generated by $\psi_{(0)}$ is $B_\Lambda$. Hence Lemma \ref{structure_lemma} part (\ref{p4}) implies that this map is a bijection to $B_\Lambda \otimes F$.

In order to prove (\ref{cabacus}) and (\ref{ccpp}), recall that the principally graded weight of a canonical basis vector $v \in V_\Lambda$ can be calculated from the crystal graph: set the weight of the highest weight element to be $0$, and let $f_i$ have degree $1$. This is a well defined grading on $B_\Lambda$ by standard results (see \cite{Hong&Kang:2000}). The principle grading on $V_\Lambda$ is obtained by letting $v(b)$ have the same weight as $b$. In particular, acting on $\psi$ by $f_i$ increases the principally graded weight of $v(\psi)$ by one.

(\ref{cabacus}): It follows from part (\ref{ccommutes}) that $\{ v(\psi) \otimes \lambda (\psi) \}$ is a basis for $V_\Lambda \otimes F$. It remains to show that the weights are correct. It is clear from the definitions that the operations $f_i$ add one to the weight of an abacus configuration, and $T_k^*$ adds $n$. As above, acting by $f_i$ increases the principally graded weight of $v(\psi)$ by $1$, without affecting $\lambda(\psi)$. Acting by $T_k^*$ adds one box to $\lambda(\psi)$ without affecting $v(\psi)$. The result follows.

(\ref{ccpp}): This follows from (\ref{cabacus}), once one notices that moving a bead one position to the right on the abacus (in such a way that it remains decreasing) always corresponds to adding $1$ to one entry of the corresponding cylindric plane partition.
\end{proof}

\subsection{Relation to the study of random cylindric plane partitions} \label{rcpp}

In \cite{Borodin:2006}, Borodin studies the expected behavior of large random cylindric plane partitions. He uses the distribution where the probability of a cylindric plane partition $\pi$ is proportional to $q^{|\pi |}$. In particular, he calculates the partition function for this system. We are now ready to present a new method to calculate this partition function. We will begin this section by explaining our method, and stating our formula. Then we will give Borodin's formula, which is quite different from ours. Next we directly show that the two formulas agree. This is not logically necessary, since the two formulas have been proven separately, but we feel it is worth including for two reasons. First, it provides a comforting verification that the results in this paper are consistent with the literature. Second, we believe that one of the most interesting consequences of our work is that it provides a new link between the study of cylindric plane partitions and the representation theory of $\asl_n$. As such, results from one area should translate to the other. We believe that directly linking the two formulas for the partition function is an important step in that direction.

The key to our method of calculating the partition function is that we can now identify the set of cylindric plane partitions on a given cylinder with a basis for $V_\Lambda \otimes F$, and that this preserves an appropriate notion of weight.  Let $E_i$ and $F_i$ be the Chevalley generators of $\g$. Recall that the principle grading of a highest weight representation $V$ of $\g$ is the $\bz_{\geq 0}$ grading induced by putting the highest weight vector in degree $0$, letting all $F_i$ have degree 1, and letting all $E_i$ have degree $-1$.  Let $V_k$ denote the degree $k$ part of $V$ in this grading. The $q$-character of $V$ is the formal sum:
\begin{equation}
\dim_q (V) = \sum_{k \geq 0} \dim (V_k) q^k.
\end{equation}
Also, recall that $F$ is the space spanned by all partitions, and define the $q^n$ character of $F$ by:

\begin{equation} \label{dimq}
\dim_{q^n} (F) := \sum_{\mbox{partitions} \hspace{1.5pt} \lambda} q^{n | \lambda |}.
\end{equation}

\begin{Theorem} \label{mypart}
The partition function for cylindric plane partitions on a cylinder whose boundary satisfies $\Lambda (\pi) = \Lambda $ (see Definition \ref{Lambdab}) is given by:
\begin{align} 
 \label{dimqs} Z := \sum_{\begin{array}{c} \mbox{cylindric partitions} \quad \pi \\ \mbox{on a given cylinder} \end{array}}  q^{| \pi |} &= \dim_q (V_\Lambda)  \dim_{q^n} (F)
 \end{align}
\end{Theorem}

\begin{proof} This follows immediately from Theorem \ref{weights} part (\ref{ccpp}).
\end{proof}

A formula for $\dim_q(V_\Lambda)$ can be found in (\cite{Kac:1990}, Proposition 10.10). The $q^n$-character of $F$ is the sum over all partitions $\lambda$ of  $q^{n | \lambda |}$, and is well known (see for example \cite{Andrews:1976}, Theorem 1.1). Substituting these results into Equation (\ref{dimqs}) we obtain an equivalent formula for $Z$:

\begin{align} \label{myZ}
 Z=
\prod_{\alpha \in \Delta_+^\vee} \left( \frac{1- q^{ \langle \Lambda + \rho, \alpha \rangle}}{1-q^{\langle \rho, \alpha \rangle}} \right)^{\mbox{mult} (\alpha)}  \prod_{k =1}^\infty \frac{1}{1-q^{kn}}.
\end{align}

In order to state Borodin's result we need some more notation.
Label the steps of the boundary of a cylindric plane partition with residues $i$ modulo $n + \ell$, as in Figure \ref{period_nl}. Define:

$\bullet$ $N= n+\ell$. This is the period of the boundary of the cylindric plane partition.

$\bullet$ For any $k \in \bz$, $k (N)$ is the smallest non-negative integer congruent to $k$ modulo $N$.

$\bullet$ $\overline{1,N}$ is the set of integers modulo $N$.

$\displaystyle \bullet A[i] = \begin{cases}
1 \quad \mbox{if the boundary is sloping up and to the right on diagonal}  \quad i \\
0 \quad \mbox{otherwise}
\end{cases}$

$\displaystyle \bullet B[i] = \begin{cases}
1 \quad \mbox{if the boundary is sloping down and to the right on diagonal}  \quad i \\
0 \quad \mbox{otherwise}
\end{cases}$

Note: We have reversed the definitions of $A[i]$ and $B[i]$ from those used by Borodin. This corresponds to reflecting the cylindric plane partition about a vertical axis, so by symmetry does not change the partition function.

\begin{Theorem} (Borodin 2006) The partition function for cylindric plane partitions is given by:
\begin{equation} \label{BorodinZ} Z  := \sum_{\begin{array}{c} \mbox{cylindric partitions} \quad \pi \\ \mbox{on a given cylinder} \end{array}}  q^{| \pi |} = 
\prod_{k \geq 1} \frac{1}{1-q^{kN}} \prod_{\begin{array}{c} i \in \overline{1,N}: A[i]=1 \\ j \in \overline{1,N}: B[j]=1 \end{array} } \frac{1}{1 - q^{(i-j)(N)+(k-1)N}}.
\end{equation}
\end{Theorem}

We will now directly show that Equations (\ref{myZ}) and (\ref{BorodinZ}) are equivalent. We already have a proof of Theorem \ref{mypart}, so some readers may wish to skip to Section \ref{wld}. 

We need the following fact about affine Lie algebras, which can be found in (\cite{Kac:1990}, Chapter 14.2). Let $\g$ be a finite dimensional simple complex Lie algebra. Let $\widehat{\g}$ be the associated untwisted affine algebra, and $\widehat{\g}'$ the corresponding derived algebra. The principle grading on $\widehat{\g}'$ is the grading determined by setting $\deg(H_i)=\deg(c)=0$, $\deg(E_i)=1$ and $\deg(F_i)= -1$, for each $i$. Let $r = \mbox{rank} (\g)$.

\begin{Lemma} (\cite{Kac:1990}, Chapter 14.2)
The dimension of the $k^{th}$ principally graded component of $\widehat{\g}'$ is $r+ s(k)$, where $s(k)$ is the number of exponents of $\g$ congruent to $k$ modulo the dual Coxeter number.
\end{Lemma}

\noindent In the case of $\asl_n$, this reduces to:

\begin{Corollary} \label{principle_dimension}
The dimension of the $k^{th}$ principally graded component of $\asl_n'$ is $n$, unless $k \equiv 0$ mod n, in which case the dimension is $n-1$.
\end{Corollary}

We start with Equation (\ref{myZ}) and manipulate it to reach Equation (\ref{BorodinZ}). First, notice that
$q^{\langle \rho, \alpha \rangle}$ is just $q$ to the principally graded weight of $\alpha$. Using this, and identifying $\Delta_+^\vee$ with $\Delta_+$ (since $\asl_n$ is self dual),  Corollary \ref{principle_dimension} allows us to simplify Equation (\ref{myZ}) as follows:

\vspace{-0.05in}

\begin{align}
Z & = \prod_{\alpha \in \Delta_+^\vee} (1- q^{ \langle \Lambda + \rho, \alpha \rangle} )^{mult ( \alpha )}
\prod_{k=1}^\infty \left( \frac{1}{1-q^k} \right)^{\begin{array}{l} \mbox{dim of the } \hspace{0.2pt} k^{th} \hspace{0.2pt} \mbox{principally} \\ \mbox{graded component of} \hspace{2pt} \asl_n' \end{array}}
\prod_{k=1}^{\infty} \frac{1}{1-q^{kn}} \label{aZeq}
\\
\label{twoparts} &=\prod_{\alpha \in \Delta_+^\vee} (1- q^{ \langle \Lambda + \rho, \alpha \rangle} )^{mult ( \alpha )}
\prod_{k=1}^\infty \left( \frac{1}{1-q^k} \right)^n.
\end{align}

We will deal with the second factor first. One can readily see that there are exactly $n$ residues $j \in\overline{1,N}$ such that $B[j]=1$ (see Figure \ref{period_nl}). Hence:

\begin{align}
\label{apart}
\prod_{k=1}^\infty \left( \frac{1}{1-q^k} \right)^n  & =  \prod_{\begin{array}{c}  j \in \{ 1, \ldots N \} , B[j]=1\\ i > j \end{array}}  \frac{1}{1-q^{i-j}}\\
\nonumber & = \prod_{k=1}^\infty \prod_{i,j \in \overline{1,N}: B[j]=1} \frac{1}{1 - q^{(i-j)(N)+(k-1)N}}
\end{align}

Now we consider the first factor of Equation (\ref{twoparts}) . Write:
\begin{equation}
\Lambda (\pi) = \sum_{i=0}^{n-1} m_i \Lambda_i.
\end{equation}
Notice that $m_i$ is the length of the $i^{th}$ downward sloping piece of the boundary of the cylindric plane partition, counting from the left over one period in Figure \ref{period_nl}, and starting with $0$ (see Definition \ref{Lambdab}). For each $i \in I$, set
$s_i := 1 + m_i.$
Let $\g_k (s)$ be the $k^{th}$ degree piece of $\asl_n'$ with the grading induced by letting $E_i$ have degree $s_i$, and $F_i$ have degree $-s_i$. This will be called the $s$-grading of $\asl_n'$. Then $q^{\langle \Lambda+\rho, \alpha \rangle}$ is just $q$ to the $s$-graded degree of $\alpha$. Again identifying $\Delta_+^\vee$ with $\Delta_+$, we see that: 

\begin{equation} \label{aroots}
 \prod_{\alpha \in \Delta_+^\vee} \left(1- q^{ \langle \Lambda + \rho, \alpha \rangle} \right)^{mult ( \alpha )} =
 \prod_{k =1}^\infty \left( 1-q^{\mbox{dim} ( \g_k (s) )} \right).
\end{equation}

\noindent $\mbox{dim} ( \g_k (s) )$ is the number of positive roots of $\asl_n$ of degree $k$ in the $s$ grading, counted with multiplicity. These are of two forms, which we must consider separately: 

\begin{enumerate} 

\item Real roots of the form $\beta+ k \delta$. Here $\beta$ is a root of $\mbox{sl}_n$, $k \geq 0$, and $k > 0$ if $\beta$ is negative. These each have multiplicity 1, and $s$ graded weight
$\langle \Lambda+ \rho, \beta \rangle + kN$ (since $\langle \Lambda + \rho, \delta \rangle = N$). 

\item Imaginary roots of the form $k \delta$. These each have multiplicity $n-1$, and $s$-graded weight $kN$.

\end{enumerate}

The imaginary roots contribute a factor of 
\begin{equation} \label{broots}
\prod_{k=1}^\infty \left( 1- q^{kN} \right)^{n-1} 
\end{equation}
to Equation (\ref{aroots}).

Next we find the contribution from the real roots. The positive roots of $\mbox{sl}_n$ are $\alpha_a + \alpha_{a+1} + \ldots + \alpha_b$ for all $1 \leq a <b \leq n-1$, and the negative roots are their negatives. A straightforward calculation shows that
\begin{equation}
\langle \rho+\Lambda , \sum_{i=a}^b \alpha_i \rangle=\left( b + \sum_{k=0}^{b-1} m_k \right)- \left(a + \sum_{k=0}^{a-1} m_k \right).\end{equation} 

One can see from Figure \ref{period_nl} that $B[a + m_0 + m_1 + \ldots + m_{a-1}]=1$ for all $a$. By a similar argument, if $\beta$ is a negative root of $\mbox{sl}_n$, then $\langle \rho+ \Lambda, \beta+\delta \rangle$ will be the difference of two integers that have $B[\cdot]=1$. It should be clear from Figure \ref{period_nl} that each pair $i,j \in \overline{1,N}$ with  $B[i]=B[j]=1$  corresponds to a root $\beta$ of $\mbox{sl}_n$ in one of these two ways. Therefore the real roots contribute a factor of
\begin{equation} \label{croots}
\prod_{k=1}^\infty \prod_{\begin{array}{c} i \in \overline{1,N}: B[i]=1 \\ j \in \overline{1,N}: B[j]=1  \\ i \neq j\end{array} } (1 - q^{(i-j)(N)+(k-1)N})
\end{equation}
to Equation (\ref{aroots}). Using Equations (\ref{broots}) and (\ref{croots}), Equations (\ref{aroots}) is equivalent to:
\begin{equation} \label{bpart}
 \prod_{\alpha \in \Delta_+^\vee} \left(1- q^{ \langle \Lambda + \rho, \alpha \rangle} \right)^{mult ( \alpha )} =
\prod_{k=1}^\infty  \left( 1- q^{kN} \right)^{n-1} \prod_{\begin{array}{c} i \in \overline{1,N}: B[i]=1 \\ j \in \overline{1,N}: B[j]=1  \\ i \neq j \end{array} } (1 - q^{(i-j)(N)+(k-1)N}).
\end{equation}
Next, substitute Equations (\ref{apart}) and (\ref{bpart}) into Equation (\ref{aZeq}) to get:
\begin{align} 
Z & = \prod_{k=1}^\infty \left( 1- q^{kN} \right)^{n-1}\prod_{\begin{array}{l} i, j \in \overline{1,N} \\ B[j]=1 \end{array}} \frac{1}{1 - q^{(i-j)(N)+(k-1)N}}
\hspace{-0.8cm}  \prod_{\begin{array}{c} i \in \overline{1,N}: B[i]=1 \\ j \in \overline{1,N}: B[j]=1  \\ i \neq j \end{array} } \hspace{-0.8cm}  \left( 1 - q^{(i-j)(N)+(k-1)N} \right) \\
& \label{almost_done} = 
\prod_{k=1}^\infty \frac{1}{1-q^{kN}}
\prod_{\begin{array}{l} i, j \in \overline{1,N} \\ B[j]=1 \end{array}} \frac{1}{1 - q^{(i-j)(N)+(k-1)N}}
  \prod_{\begin{array}{c} i \in \overline{1,N}: B[i]=1 \\ j \in \overline{1,N}: B[j]=1  \end{array} } \hspace{-0.8cm}  \left( 1 - q^{(i-j)(N)+(k-1)N} \right).
\end{align}
Equation (\ref{almost_done}) follows because there are exactly $n$ residues $i \in \overline{1,N}$ that have $B[i]=1$. This then simplifies to Equation (\ref{BorodinZ}), using the fact that $A[i]=1$ if and only if $B[i] \neq 1$. Hence we have directly shown that our equation for the partition function agrees with Borodin's.

\subsection{A rank-level duality result} \label{wld}
We now know that the set of cylindric plane partitions with a given boundary parameterizes a level $\ell$ crystal for $\asl_n$. We can reflect Figure \ref{period_nl} about a vertical axis to see that they just as easily parametrize a level $n$ crystal for $\asl_\ell$. Equation (\ref{dimqs}) then gives two ``dual" ways to calculate the partition function $Z$. Hence we get an identity involving the $q$-characters of the two corresponding representations. We state this precisely below. We then show that our result is in fact a variant of a result due to Frenkel \cite{IFrenkel:1982}.   

For a residue $i$ modulo $n$, let  $\Lambda_i^{(n)}$ denote the $i^{th}$ fundamental weight of $\asl_n$.
For any level $\ell$ dominant integral weight $\Lambda = \sum_{i=0}^{n-1} c_i \Lambda_i^{(n)}$ of $\asl_n$, define a corresponding level $n$ dominant integral weight of $\asl_\ell$ by:
\begin{equation} \label{Lambda'}
\Lambda' := \sum_{i=0}^{n-1} \Lambda^{(\ell)}_{c_i+ c_{i+1}+ \cdots + c_{n-1}}.
\end{equation}
This should be understood using Figure \ref{period_nl}: If $\Lambda$ is defined as in the caption, then $\Lambda'$ is defined in the same way, but after first reflecting in the vertical axis, interchanging diagonal $\pi_k$ and diagonal $c_k$.

\begin{Theorem} \label{wldtheorem}
Let $\Lambda$ be a level $\ell$ dominant integral weight of $\asl_n$, with $n, \ell \geq 2$. Let $\Lambda'$ be the corresponding dominant integral weight of $\asl_\ell$ defined by Equation (\ref{Lambda'}). Then:
\begin{equation} \label{wldequation}
\dim_q (V_\Lambda) \prod_{k=0}^\infty \frac{1}{1-q^{nk}} = \dim_q (V_{\Lambda'}) \prod_{k=0}^\infty \frac{1}{1-q^{\ell k}}.
\end{equation}
\end{Theorem}

\begin{proof}
Consider the system $S$ consisting of all cylindric plane partitions with a fixed boundary, where the probability of $\pi$ is proportional to $q^{|\pi|}$, and the boundary is chosen so that $\Lambda(\pi)= \Lambda$ for any $\pi \in S$. By Theorem \ref{mypart}, the partition function for $S$ is given by:
\begin{equation} \label{aa1}
Z := \sum_{\pi \in S}  q^{|\pi|} = \dim_q (V_\Lambda) \dim_{q^n} (F). 
\end{equation}
Next consider the system $S'$ consisting of all cylindric plane partitions with a fixed boundary, where the probability of $\pi$ is proportional to $q^{|\pi|}$, and the boundary is chosen so that $\Lambda(\pi)= \Lambda'$ for any $\pi \in S'$. Once again, by Theorem \ref{mypart}, the partition function for $S'$ is given by:
\begin{equation} \label{aa2}
Z' := \sum_{\pi \in S'}  q^{|\pi|} = \dim_q (V_{\Lambda'}) \dim_{q^\ell} (F). 
\end{equation}
As discussed in the caption of Figure \ref{period_nl}, we can reflect an element of $S$ about a vertical axis to get an element of $S'$. This is a weight preserving bijection between $S$ and $S'$, so the two partition functions $Z$ and $Z'$ must be equal. The result then follows from Equations (\ref{aa1}) and (\ref{aa2}), substituting $\displaystyle \sum_{r \geq 1} \frac{1}{1-s^r}$ for $\dim_s(F)$ on both sides.
\end{proof}
For completeness, we also include the corresponding result when $\ell=1$:

\begin{Theorem} \label{n2}
Let $\Lambda$ be a level $1$ dominant integral weight of $\asl_n$, with $n \geq 2$. Then:
\begin{equation}
\label{n22} \dim_q (V_\Lambda) \prod_{k=0}^\infty \frac{1}{1-q^{nk}} = \prod_{k=0}^\infty \frac{1}{1-q^k}.
\end{equation}
\end{Theorem}

\begin{proof}
By Theorems \ref{mypart} and \ref{acppb}, the left side is $\sum q^{|\psi|}$, where the sum is over all $1$-strand (descending) abacus configurations $\psi$ with a given compactification $\psi_{(0)}$. Such a configuration is just one row of beads. We can simultaneously shift each $\psi$ until $\psi_{(0)}$ has its last black bead in position $-1/2$, without changing any $|\psi|$. As in Section \ref{lotsastuff}, there is a bijection between partitions and rows of beads whose compactification is this $\psi_{(0)}$, and one can easily see that this bijection preserves the weights. Hence the left side of Equation (\ref{n22}) is equal to $\sum_{\mbox{partitions } \lambda} q^{|\lambda|},$ which is well known to be given by the right side. \end{proof}

In \cite{IFrenkel:1982}, Frenkel proves that the character of an irreducible integral level $\ell$ representation of $\hgl_n$ is equal to the character of a corresponding irreducible integral level $n$ representation of $\hgl_\ell$. His proof is a direct argument using the Kac-Weyl character formula. We now show that Frenkel's duality is equivalent to Theorem \ref{wldtheorem}, thereby obtaining a new proof of Frenkel's result. First we set up some notation. Our definitions agree with those in \cite{IFrenkel:1982}, although we have reworded some parts.

We use the following realization of $\hgl_n$:
\begin{equation}
\hgl_n = \mbox{gl}_n \otimes \bc[t, t^{-1}] + \bc c + \bc d.
\end{equation}

\noindent Here $c$ is central and $d$ acts by the derivation $t \frac{d}{dt}$. The bracket on the remaining part is defined by
\begin{equation}
[x \otimes t^p, y \otimes t^q] = [x,y] \otimes t^{p+q} + p \delta_{p+q,0} \langle x, y \rangle c,
\end{equation}
where $\langle x, y \rangle= \mbox{Tr} ( \mbox{ad}(x) \mbox{ad}(y)).$

Fix a dominant integral weight $\Lambda$ for $\asl_n$. We say an irreducible representation $W_\Lambda$ of $\hgl_n$ is a highest weight representation of highest weight $\Lambda$ if there is some $w_\Lambda \in W_\Lambda$ such that:
\begin{enumerate}

\item{  $w_\Lambda$ is a highest weight vector of $W_\Lambda$ as an $\asl_n$ representation, and has weight $\Lambda$.}

\item{$\mbox{Id} \otimes t^n \cdot w_\Lambda=0$ for all $n >0$.}

\end{enumerate}

The principle grading on such a $W_\Lambda$ can be defined as follows: Put $w_\Lambda$ in degree zero. For each $0 \leq i \leq n-1$, let $E_i$ have degree $-1$ and $F_i$ have degree 1. Finally, let $c$ and $d$ have degree $0$, and $\mbox{Id} \otimes t^k$ have degree $-nk$. Denote by $W_{\Lambda,k}$ the degree $k$ weight space in this grading. The $q$-character of $W_\Lambda$ is the formal sum:
\begin{equation} \label{qqchar}
\dim_q(W_\Lambda)= \sum_{k=0}^\infty \dim(W_{\Lambda,k})q^k.
\end{equation}

\begin{Comment} \label{abcd}
It is a straightforward exercise to see that this definition agrees with that given in \cite{IFrenkel:1982}, and that the restriction of this grading to the $\asl_n$ representation $\asl_n \cdot w_\Lambda$ agrees with the principle grading defined in Section \ref{rcpp}. 
One should also note that for any given $\Lambda$ there are many possibilities for $W_\Lambda$, since  $\mbox{Id} \otimes 1$ can act by any scalar. However, these all have the same $q$-character, so the $q$-character is well defined as a function of $\Lambda$.
\end{Comment}

The following result is due to Frenkel  (\cite{IFrenkel:1982}, Theorem 2.3):

\begin{Theorem} (Frenkel 1982) \label{IF:82}
Fix a level $\ell$ dominant integral  weight $\Lambda$ for $\asl_n$, and let $\Lambda'$ be the corresponding level $n$ dominant integral weight of $\asl_\ell$ given by Equation (\ref{Lambda'}). Let $W_\Lambda$ be an irreducible representation of $\hgl_n$ of highest weight $\Lambda$, and $W_{\Lambda'}$ be an irreducible representation of $\hgl_\ell$ of highest weight $\Lambda'$. Then:
\begin{equation}
\dim_q(W_\Lambda)= \dim_q(W_{\Lambda'}).
\end{equation}
\end{Theorem}

The following lemma will give us an alternative proof of Theorem \ref{IF:82}:

\begin{Lemma} \label{aformula}
Fix a dominant integral weight $\Lambda$ for $\asl_n$. Let $V_\Lambda$ be the irreducible $\asl_n$ representation of highest weight $\Lambda$, and $W_\Lambda$ be an irreducible highest weight representation of $\hgl_n$ of highest weight $\Lambda$. Then:
\begin{equation} 
\dim_q(W_\Lambda)= \dim_q(V_\Lambda) \prod_{k=1}^\infty \frac{1}{1-q^k}.
\end{equation}
\end{Lemma}

\begin{proof}

Let $\mathcal{W}$ be the $\bc$ span of $c$ and $\mbox{Id} \otimes t^p$ for $p \in \bz$. Let $\mathcal{S}$ be the $\bc$ space of $c$, $d$, and $x \otimes t^p$ for $x \in \mbox{sl}_n$ and $p \in \bz$.
One can check that $\mathcal{W}$ forms a copy of the infinite dimensional Weyl algebra and $\mathcal{S}$ forms a copy of $\asl_n$. Furthermore, $[\mathcal{W},\mathcal{S}]=0$, $\hgl_n = \mathcal{W} + \mathcal{S}$, and $\mathcal{W} \cap \mathcal{S} = \bc c$. It follows that
\begin{equation} \label{stufff}
W_\Lambda= V \boxtimes F,
\end{equation}
where $V$ is an irreducible representation of $\asl_n$ and $F$ is an irreducible representation of $\mathcal{W}$. By Comment \ref{abcd}, we must have $V=V_\Lambda$. There are many irreducible highest weight representations of $\mathcal{W}$ (parametrized by the actions of $c$ and $\Id \otimes t^0$), but each has the same $q$-character with respect to the grading where $c$ has degree $0$, $\Id \otimes t^k$ has weight $-n k $, and the highest weight space is placed in degree $0$. This is given by:
\begin{equation}
\dim_q(F) = \prod_{k=1}^\infty \frac{1}{1-q^{nk}}.
\end{equation}

The Lemma follows by taking the $q$-character of each side of Equation (\ref{stufff}).
\end{proof}

Theorem \ref{IF:82} follows by applying Lemma \ref{aformula} to each side of Equation (\ref{wldequation}). \qed

\section{Relation to the Kyoto path model} \label{like_kashiwara}

We now present an explicit bijection between the set of tight descending abacus configurations with a given compactification $\psi_{(0)}$, and the Kyoto path model for an integrable irreducible level $\ell$ representation of $\asl_n$, and show that the images of our operators $e_i$ and $f_i$ are the crystal operators $e_i$ and $f_i$. This proof does not assume that the abacus model is an $\asl_n$ crystal. Since the Kyoto path model is an $\asl_n$ crystal, we get a new proof of Theorem \ref{hipart}, which works for $\asl_2$ (where our previous proof failed). Since Lemma \ref{structure_lemma} parts (\ref{p2}) and (\ref{p4}) also hold in the case of $\asl_2$, the set of all descending abacus configurations is still an $\asl_2$ crystal. Hence the results in Section \ref{cylindric_partitions_and_abacus} still hold in the case $n=2$.

We refer the reader to (\cite{Hong&Kang:2000}, Chapter 10) for the details of the Kyoto path model. We will only use the path models corresponding  to one family $B_\ell$ of perfect crystals for $\asl_n$, one of each level $\ell > 0$. These can be found in (\cite{KKMMNN2} Theorem 1.2.2), where they are called $B^{1,\ell}$. Here $B_\ell$ consists of semi standard fillings of the partition $( \ell )$ with $0.5, 1.5, \ldots (2n-1)/2$. Below is an example for $n=3$ and $\ell=4$:

\begin{center}
\begin{picture}(4,1)
\put(0,0){\line(1,0){4}}
\put(0,1){\line(1,0){4}}
\put(0,0){\line(0,1){1}}
\put(1,0){\line(0,1){1}}
\put(2,0){\line(0,1){1}}
\put(3,0){\line(0,1){1}}
\put(4,0){\line(0,1){1}} 
\put(0.1, 0.32){\small{0.5}}
\put(1.1, 0.32){\small{1.5}}
\put(2.1, 0.32){\small{2.5}}
\put(3.1, 0.32){\small{2.5}}
\end{picture}

\end{center}

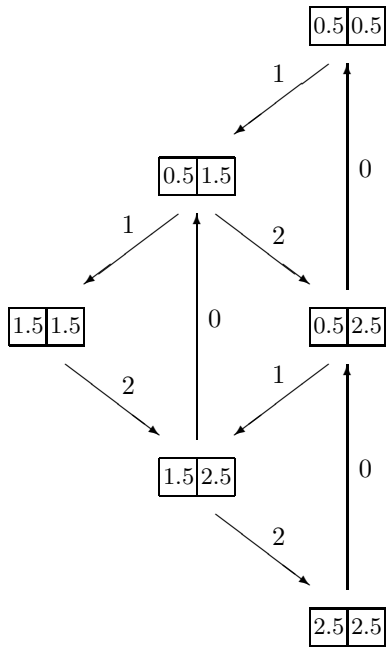
\begin{figure}
\begin{picture}(30,18)

\put(17.5, 15.5){\vector(-4,-3){2.5}}
\put(13.5, 11.5){\vector(-4,-3){2.5}}
\put(17.5, 7.5){\vector(-4,-3){2.5}}
\put(16,15){1}
\put(12,11){1}
\put(16,7){1}
\put(14.5, 11.5){\vector(4,-3){2.5}}
\put(10.5, 7.5){\vector(4,-3){2.5}}
\put(14.5, 3.5){\vector(4,-3){2.5}}
\put(16,10.7){2}
\put(12,6.7){2}
\put(16,2.7){2}
\put(18,9.5){\vector(0,1){6}}
\put(18,1.5){\vector(0,1){6}}
\put(14,5.5){\vector(0,1){6}}
\put(18.3,12.5){0}
\put(18.3,4.5){0}
\put(14.3,8.5){0}

\put(17,16){\begin{picture}(2,1)
\put(0,0){\line(1,0){2}}
\put(0,1){\line(1,0){2}}
\put(0,0){\line(0,1){1}}
\put(1,0){\line(0,1){1}}
\put(2,0){\line(0,1){1}} 
\put(0.1, 0.32){\small{0.5}}
\put(1.1, 0.32){\small{0.5}}
\end{picture}}

\put(13,12){\begin{picture}(2,1)
\put(0,0){\line(1,0){2}}
\put(0,1){\line(1,0){2}}
\put(0,0){\line(0,1){1}}
\put(1,0){\line(0,1){1}}
\put(2,0){\line(0,1){1}} 
\put(0.1, 0.32){\small{0.5}}
\put(1.1, 0.32){\small{1.5}}
\end{picture}}

\put(9,8){\begin{picture}(2,1)
\put(0,0){\line(1,0){2}}
\put(0,1){\line(1,0){2}}
\put(0,0){\line(0,1){1}}
\put(1,0){\line(0,1){1}}
\put(2,0){\line(0,1){1}} 
\put(0.1, 0.32){\small{1.5}}
\put(1.1, 0.32){\small{1.5}}
\end{picture}}

\put(17,8){\begin{picture}(2,1)
\put(0,0){\line(1,0){2}}
\put(0,1){\line(1,0){2}}
\put(0,0){\line(0,1){1}}
\put(1,0){\line(0,1){1}}
\put(2,0){\line(0,1){1}} 
\put(0.1, 0.32){\small{0.5}}
\put(1.1, 0.32){\small{2.5}}
\end{picture}}

\put(13,4){\begin{picture}(2,1)
\put(0,0){\line(1,0){2}}
\put(0,1){\line(1,0){2}}
\put(0,0){\line(0,1){1}}
\put(1,0){\line(0,1){1}}
\put(2,0){\line(0,1){1}} 
\put(0.1, 0.32){\small{1.5}}
\put(1.1, 0.32){\small{2.5}}
\end{picture}}

\put(17,0){\begin{picture}(2,1)
\put(0,0){\line(1,0){2}}
\put(0,1){\line(1,0){2}}
\put(0,0){\line(0,1){1}}
\put(1,0){\line(0,1){1}}
\put(2,0){\line(0,1){1}} 
\put(0.1, 0.32){\small{2.5}}
\put(1.1, 0.32){\small{2.5}}
\end{picture}}

\end{picture}
\caption{A perfect crystal of level 2 for $\asl_3$. The edges show the actions of $f_0, f_1,f_2$. The $e_i$ just follow the arrows backwards. If there is no $i$ arrow coming out of a vertex, $f_i$ sends that element to $0$. \label{perfect_crystal}}
\end{figure}

The operator $f_i$ changes the rightmost $i - 1/2$ to $i + 1/2$, if possible, and sends the element to 0 if there is no  $i - 1/2$. $f_0$ changes one $n - 1/2$ to $1/2$, if possible, then shuffles that element to the front. If there is no $n-1/2$ then $f_0$ sends that element to $0$. $e_i$ inverts $f_i$ if possible, and send the element to $0$ otherwise. This should be clear from Figure \ref{perfect_crystal}. Notice that we are using half integers instead of integers. The usual conventions are obtained by adding $1/2$ to everything.

Given a level $\ell$ perfect crystal $B_\ell$ for an affine algebra $\g$, and a level $\ell$ integral weight $\lambda$ of $\g$, 
Kashiwara et. al. define a {\it ground state path} of weight $\lambda$ to be a sequence $${\bf p}_\lambda =
\cdots \otimes b_{3}  \otimes  b_2 \otimes b_1 \in \cdots \otimes B_\ell \otimes B_\ell \otimes B_\ell$$ such that:
\begin{enumerate}
\item  $\varphi ( b_1) = \lambda$
\item For each $k \geq 1$, $\varepsilon ( b_k) = \varphi (b_{k+1})$.
\end{enumerate}
They show that, once $B_\ell$ is fixed, there is a unique ground state path of weight $\lambda$ for any level $\ell$ dominant integral weight $\lambda$. 

The irreducible crystal $P(\lambda)$ of highest weight $\lambda$ is the set of all all paths in $\cdots \otimes B_\ell \otimes B_\ell \otimes B_\ell$ which differ from the ground state path in only finitely many places. The crystal operators $e_i$ and $f_i$ are obtained from $B_\ell$ using the tensor product rule.
We use the notation $P_{n, \ell} (\Lambda)$ to denote the Kyoto path model for the $\asl_n$ crystal $B_\Lambda$ using the perfect crystal $B_\ell$ defined above (where $\ell$ is the level of $\Lambda$). 

Given an abacus configuration $\psi$, we define a sequence $\overline{\psi_\bigdot^k}$ of elements in $B_\ell$ (where $\ell$ is the number of rows of the abacus) by letting $\overline{\psi_\bigdot^k}$ be the unique semi-standard filling of $(\ell)$ with $\{ \psi_0^k(n), \ldots, \psi_{\ell-1}^k(n) \}$. Recall that $\psi_j^k$ is the position of the $k^{th}$ black bead from the right on row $j$ of $\psi$, and, as in Section \ref{rcpp}, $\psi_j^k(n)$ is the unique number $0 < \psi_j^k(n) < n$ that is congruent to $\psi_j^k$ modulo $n$.

\begin{Theorem} \label{Like_Kashiwara} Consider an $\ell$ row abacus colored with $c_0, \ldots c_{n-1}$ as in figure \ref{abacuscrystal}. Let $B$ be the set of tight, descending abacus configurations with compactification $\psi_{(0)}$. The map $J: B \rightarrow P_{n,\ell} (\Lambda (\psi_{(0)}))$ given by $\psi \rightarrow \cdots \otimes \overline{\psi_\bigdot^3} \otimes \overline{\psi_\bigdot^2} \otimes \overline{\psi_\bigdot^1}$ is a bijection. Furthermore, the images of our operators $e_i$ and $f_i$ are the crystal operators on $P_{n,\ell} (\Lambda(\psi_{(0)}))$.
\end{Theorem}

\begin{Comment}
This theorem really says $J$ is a crystal isomorphism. It gives a new proof of Theorem \ref{hipart}, which works even in the case $n=2$.
\end{Comment}

\begin{proof}

First note that the image of $\psi_{(0)}$ will satisfy the conditions to be a ground state path of weight $\Lambda(\psi_{(0)})$. For any $\psi \in B$, $\psi_j^k={\psi_{(0)}}_j^k$ for all but finitely many $0 \leq j \leq \ell-1$ and $k \geq 1$. Hence $J(\psi)$ will always lie in $P_{n,\ell} (\Lambda (\psi_{(0)}))$.

Use the string of brackets $S_i$ from Lemma \ref{gglemma} to calculate $f_i$ acting on $B$, and use the string of brackets $S_i'$ from Corollary \ref{mycrystaldef} to calculate $f_i$  acting on  $P_{n,\ell}(\Lambda(\psi)).$ One can easily see that $S_i'(J(\psi))=S_i(\psi)$. It follows from the two rules that the image of our operator $f_i$ is the crystal operator $f_i$ on $P_{n,\ell} (\Lambda(\psi_{(0)}))$. The proof for $e_i$ is similar.

Our operators $e_i$ and $f_i$ act transitively on $B$ by Lemma \ref{structure_lemma} part (\ref{p3b}), and the crystal operators $e_i$ and $f_i$ act transitively on  $P_{n,\ell} (\Lambda(\psi_{(0)}))$ since it is an irreducible crystal. $J$ preserves these operators, so it must be a bijection.
\end{proof}

\section{Questions} \label{questions}

We finish with a short discussion of two questions which we feel are natural next steps for the work presented in this paper.

\begin{Question}
Can one lift our crystal structures to get representations of $U_q(\asl_n)$? In particular, do cylindric plane partitions parametrize a basis for a representations of $U_q(\asl_n)$ in any natural way?
\end{Question}

In \cite{KMPY:1996}, Kashiwara, Miwa, Petersen and Yung study a $q$-deformed Fock space. This space has an action of $U_q(\asl_n)$, and a commuting action by the Bosons. In our picture, the space spanned by cylindric plane partitions has a $\asl_n$ crystal structure and a commuting action of $\gli$. Using the Boson-Fermion correspondence, one can realize $\gli$ inside a completion of the Bosonic algebra. We hope there is a natural action of $U_q(\asl_n)$ on the space spanned by cylindric plane partitions, coming from an embedding of this space into the $q$-deformed Fock space. 

This should have interesting consequences because of our rank-level duality: The space spanned by cylindric plane partitions on a given cylinder would simultaneously carry actions of $U_q(\asl_n)$, $U_q(\asl_\ell)$, and two copies of $\gli$. One copy of $\gli$ would commute with the $U_q(\asl_n)$ action and the other would commute with the $U_q(\asl_\ell)$ action. It would be interesting to see how the other actions interact. 

There are also several papers of Denis Uglov which may shed some light on this question (see \cite{Uglov:1999}, and references therein). He studies an action of $\asl_n$ on the space spanned by $\ell$-tuples of charged partitions. The crystal structure coming from his action is very similar to ours, and we believe it can be made to agree exactly by changing some conventions. He also studies commuting actions of $\asl_\ell$ and the Heisenberg algebra. It would be nice to fit these into our picture as well. However, his action does not appear to preserve the space spanned by descending abacus configurations, so it will not descend directly to the desired action on the space spanned by cylindric plane partitions.

\begin{Question}
Our Theorem \ref{mypart} gives a new way to calculate the partition function for a system of random cylindric plane partitions, using the representation theory of $\asl_n$. Can one use similar methods to do other calculations on such a statistical system? In particular, does this give a new way to calculate the correlation functions found by Borodin in \cite{Borodin:2006}? Alternatively, can one use results about random cylindric plane partitions to get interesting results about the representation theory of $\asl_n$?

\end{Question}

We have not pursued this line of inquiry in any depth, but it seems we have a new link between cylindric plane partitions and $\asl_n$ representation theory. One can certainly hope that results from one area will transfer in some way to the other.

\bibliographystyle{plain}
\bibliography{mybib}

\end{document}